\documentclass[11pt,letterpaper] {amsart}

\usepackage{graphicx,amsthm}
\usepackage{verbatim,amsmath,amscd,amssymb,color}
\usepackage[mathscr]{eucal}
\usepackage{xcolor}
\usepackage{verbatim}
\usepackage{tikz}
\usepackage{mathrsfs} 

\colorlet{mycol}{black}

\usepackage{txfonts}
\definecolor{navyblue}{RGB}{20,20,140}

\usepackage{epic,eepic,multicol}
\numberwithin{equation}{section}
\input xy
\xyoption{all}

\def\m@th{\mathsurround=0pt}
\def\fsquare(#1,#2){
\hbox{\vrule$\hskip-0.4pt\vcenter to #1{\normalbaselines\m@th
\hrule\vfil\hbox to #1{\hfill$\scriptstyle #2$\hfill}\vfil\hrule}$\hskip-0.4pt
\vrule}}

\newcommand{\cA}{{\mathcal A}}

\newcommand{\cD}{{\mathcal D}}

\newcommand{\cH}{{\mathcal H}}
\newcommand{\cK}{{\mathcal K}}

\newcommand{\cS}{{\mathcal S}}
\newcommand{\cT}{{\mathcal T}}

\newcommand{\cP}{{\mathcal P}}

\newcommand{\fra}{\mathfrak a}

\newcommand{\frg}{\mathfrak g}

\newcommand{\frn}{\mathfrak n}

\newcommand{\frt}{\mathfrak t}

\newcommand{\bbA}{\mathbb A}
\newcommand{\bbB}{\mathbb B}
\newcommand{\bbC}{\mathbb C}

\newcommand{\bbQ}{\mathbb Q}

\newcommand{\bbZ}{\mathbb Z}

\newcommand{\scC}{{\mathscr C}}

\newcommand{\ch}{\mathrm{ch}}

\newcommand{\Hom}{\mathrm{Hom}}

\newcommand{\id}{\mathrm{id}}

\newcommand{\wt}{\mathrm{wt}}

\newcommand{\lbr}{\begin{bmatrix}}
\newcommand{\rbr}{\end{bmatrix}}
\newcommand{\cd}{commutative diagram }

\def\ge{\frg}

\def\al{\alpha}
\def\aq{\cA_q(\frn)}
\def\tilaq{\widetilde{\cA_q}(\frn)}
\def\bfk{{\bf k}}
\def\ci(#1,#2){c_{#1}^{(#2)}}
\def\Ci(#1,#2){C_{#1}^{(#2)}}
\def\mpp(#1,#2,#3){#1^{(#2)}_{#3}}
\def\bCi(#1,#2){\ovl C_{#1}^{(#2)}}
\def\ch(#1,#2){c_{#2,#1}^{-h_{#1}}}
\def\cc(#1,#2){c_{#2,#1}}

\def\bfi{{\mathbf i}}

\def\bfh{{\mathbf h}}

\def\beneme{\begin{enumerate}}
\def\beq{\begin{equation}}
\def\beqn{\begin{eqnarray}}
\def\beqnn{\begin{eqnarray*}}

\def\bbra#1,#2,#3{\left\{\begin{array}{c}\hspace{-5pt}
#1;#2\\ \hspace{-5pt}#3\end{array}\hspace{-5pt}\right\}}
\def\cd{\cdots}
\def\ci(#1,#2){c_{#1}^{(#2)}}

\def\del{\delta}
\def\Del{\Delta}

\def\eit{\tilde{e}_i}
\def\Eit{\widetilde{E}_i}
\def\eneme{\end{enumerate}}

\def\eeq{\end{equation}}
\def\eeqn{\end{eqnarray}}
\def\eeqnn{\end{eqnarray*}}
\def\fit{\tilde{f}_i}
\def\Fit{\widetilde{F}_i}

\def\gau#1,#2{\left[\begin{array}{c}\hspace{-5pt}#1\\
\hspace{-5pt}#2\end{array}\hspace{-5pt}\right]}

\def\HOM{{\rm H\textsc{om}}}

\def\ify{\infty}
\def\id{{\rm id}}

\def\lan{\langle}

\def\max{{\rm max}}
\def\lm{\lambda}
\def\Lm{\Lambda}
\def\mapright#1{\smash{\mathop{\longrightarrow}\limits^{#1}}}

\def\nd{\noindent}
\def\nn{\nonumber}

\def\ot{\otimes}

\def\ovl{\overline}

\def\qq{\qquad}
\def\q{\quad}
\def\qed{\hfill\framebox[2mm]{}}

\def\ran{\rangle}
\def\rgmod{\hbox{$R$-\hbox{gmod}}}

\def\ssl{\mathfrak sl}

\def\til{\tilde}
\def\tilrgmod{\widetilde{R}\hbox{-gmod}}

\def\tt{\frt}

\def\uq{U_q(\ge)}

\def\uqm{U^-_q(\ge)}
\def\uqp{U^+_q(\ge)}

\def\TY(#1,#2,#3){#1^{(#2)}_{#3}}

\def\vep{\varepsilon}
\def\vp{\varphi}

\def\xxi(#1,#2,#3){\displaystyle {}^{#1}\Xi^{(#2)}_{#3}}

\def\wtil{\widetilde}
\def\what{\widehat}

\def\m@th{\mathsurround=0pt}
\def\fsquare(#1,#2){
\hbox{\vrule$\hskip-0.4pt\vcenter to #1{\normalbaselines\m@th
\hrule\vfil\hbox to #1{\hfill$\scriptstyle #2$\hfill}\vfil\hrule}$\hskip-0.4pt
\vrule}}

\theoremstyle{definition}
\newtheorem{df}{Definition}[section]
\newtheorem{thm}[df]{Theorem}
\newtheorem{pro}[df]{Proposition}
\newtheorem{lem}[df]{Lemma}
\newtheorem{ex}[df]{Example}

\newtheorem{rem}[df]{Remark}

\title[Categorified Crystal Structure on Localized Quantum Coordinate Rings]
{ Categorified Crystal Structure on Localized Quantum Coordinate Rings}

\author{ T\textsc{oshiki} N\textsc{akashima}}\thanks{
{Division of Mathematics, 
Sophia University, Kioicho 7-1, Chiyoda-ku, Tokyo 102-8554,
Japan}\\
Email:\texttt{toshiki@sophia.ac.jp},\,\,
T.N is supported in part by
JSPS Grants in Aid for Scientific Research $\#$20K03564, }

\date{}

\begin{document}
\maketitle
\begin{abstract}
For the quiver Hecke algebra $R$ associated with a simple Lie algebra, 
let $\rgmod$ be the category of finite-dimensional 
graded $R$-modules. It is well-known that it 
categorifies the unipotent quantum coordinate ring $\aq$. The 
localization of  $\rgmod$ has been defined in \cite{Loc}, 
which will be denoted by $\tilrgmod$. 
Its Grothendieck ring $\cK(\tilrgmod)$ defines the localized (unipotent) 
quantum coordinate ring
 $\wtil\aq$. We shall give a certain crystal structure  
 on the localized quantum coordinate ring by 
regarding the set of self-dual simple objects $\bbB(\tilrgmod)$ in $\tilrgmod$.
We also give the isomorphism of crystals from $\bbB(\tilrgmod)$ to the cellular crystal
$\bbB_\bfi$ for an arbitrary reduced word $\bfi$ of the longest Weyl group element. 
This result can be seen as a localized version 
for the categorification of the crystal
$B(\ify)$ by Lauda-Vazirani \cite{L-V} since the crystal 
$B(\ify)$ is realized in the cellular crystal $\bbB_\bfi$.

\end{abstract}

%%%%%%%%%%%%%%%%%%%%%%%%%%%%%%%%%%%%%%%%
\section{Introduction}

The quiver Hecke algebra has been introduced by 
Khovanov-Lauda (\cite{K-L,K-L2}) and Rouquier (\cite{Rou}) independently, which is a family of $\bbZ$-graded ${\bf k}$-algebras
$R(\beta)$ defined for $\beta=\sum_im_i\al_i\in Q_+$, where $Q_+$ is the 
positive root lattice and ${\bf k}$ is a field.
For a symmetrizable Kac-Moody Lie algebra $\ge$, 
the remarkable property on the algebra $R:=\bigoplus_{\beta\in Q_+}R(\beta)$ is that 
it categorifies the nilpotent half of the quantum algebra $\uqm$ and 
the unipotent quantum coordinate ring $\aq$, that is, 
there exist isomorphisms of algebras:
\[
\cK(R\hbox{\rm-proj})\cong U^-_q(\ge)_{\bbZ[q,q^{-1}]}\qq\qq\cK(\rgmod)\cong 
\cA(\mathfrak n)_{\bbZ[q,q^{-1}]},
\]
where $\cK(R\hbox{\rm-proj})$ (resp. $\cK(\rgmod)$) is the Grothendieck ring of 
$R\hbox{\rm-proj}$ (resp. $\rgmod$) the category of graded projective(resp. finite-dimensional) $R$-modules, which is equipped with a multiplication by
the convolution product on modules. In particular, 
in \cite{Rou2,VV} for a symmetric $\ge$  this isomorphism yields one to one correspondence between 
 the set of indecomposable projective $R$-modules 
(resp.  finite-dimensional simple $R$-modules) and 
the set of the lower global bases/canonical bases 
(resp. upper global bases/dual canonical bases) (\cite{K1},\cite{Lus}). 
Those were first introduced in \cite{Lus} for the ADE cases.

The theory of Kashiwara's crystal bases \cite{K1,K3} involved
various applications to many areas, such as, mathematical physics
\cite{KMN1,KMN2}, combinatorics \cite{KN,Lit,N0}, cluster algebras \cite{FZ1}, 
modular representations and LLT conjecture
\cite{Ariki,LLT}, etc. 
In particular, in \cite{Ariki} Ariki has just open the door to
 the method of categorification, which has become the standard 
in representation theory nowadays.
The categorification as above also played a significant role
in realizing the cluster algebra on  the unipotent quantum coordinate rings, which is called a 
monoidal categorification of cluster algebras.

Theory of crystal bases would be integrated into {\it crystals},
which  is a certain combinatorial object holding similar properties of crystal bases and 
consists of 6-tuple $(B,\wt, \{\vep_i\},\{\vp_i\}, \{\eit\},\{\fit\})_{i\in I}$
with a set $B$, the maps $\eit,\fit:B\sqcup\{0\}\to B\sqcup\{0\}$ and 
the functions $\vep_i,\vp_i:B\to \bbZ\sqcup\{-\ify\}$ and $\wt:B\to P$
 satisfy certain conditions (see Definition \ref{cryst}).
Lauda and Vazirani \cite{L-V} introduced the crystal structure on the set of finite dimensional simple $R$-modules and obtained the 
outstanding result that this crystal is isomorphic to the crystal $B(\ify)$ 
of the nilpotent half of the quantum algebra $\uqm$, that is, it can be said that
they succeeded in categorifying the crystal $B(\ify)$ in terms of the $R$-modules.

Theory of crystal base gave rise to another important object "geometric crystal"\cite{BK}, 
which is a kind of geometric lifting of crystal base  and 
is defined by replacing the components in the 6-tuple of 
crystal such as  
the set $B$ by a variety $X$, the integer-function $\vep_i$ 
by rational function $\til\vep_i$ and 
the map $\eit$ (or $\fit$) by some unital rational $\bbC^\times$-
action $e_i^c$ on $X$ ($c\in \bbC^\times$) satisfying the specific conditions.
One of the most remarkable feature of geometric crystals 
is that some geometric crystals possess "positive structure" and it yields 
"tropicalization functor", by which one has a Langlands dual crystal.
Indeed, in \cite{N3} we considered the geometric crystals on the Schubert cell
$X_w$ associated with a Weyl group element $w$ and applying the 
tropicalization/ultra-discretization functor 
we obtained a Langlands dual 
crystal $\bbB_{\bfi}^\vee$, where $\bfi=i_1\cd i_k$ is a reduced word for $w$ and this crystal called a "cellular crystal" associated with a reduced word $\bfi$, which is defined as follows: 
Let  $B_i=\{(n)_i\mid n\in\bbZ\}$
($i\in I$) be the crystal as in Example \ref{ex-cry} and then define
$\bbB_{i_1i_2\cd i_k}:=B_{i_1}\ot B_{i_2}\ot\cd\ot B_{i_k}$. 
Indeed, since the tensor product of crystals mean a direct product as sets
(see Proposition \ref{tensor}),
one finds that $\bbB_{i_1i_2\cd i_k}$ can be identified with $\bbZ^k$ as a set.
Thus, its crystal structure is rather simpler than the one of $B(\ify)$
the crystal base of the nilpotent half of the quantum algebra
$\uqm$, which is realized in some cellular crystal $\bbB_{i_1i_2\cd i_k}$ (see \cite{N2,N-Z})
if $\ge$ is of finite type.
The crystal $B_i$ and $\bbB_\bfi$ are sometimes called  virtual crystals, 
which means that they are not obtained by crystallization process from 
any representation-theoretical object, such as, $\uqm$ or irreducible 
highest weight modules $V(\lm)$ but defined as abstract combinatorial objects.
Nevertheless, the tropicalization process admits a certain representation-theoretical
background on such "virtual crystals".

In \cite{Loc} (see also \cite{KKK}), 
Kashiwara et al have invented the "localization" functor on 
graded monoidal categories. 
For a graded monoidal category $(\cT=\oplus_{\lm\in \Lm}\cT_\lm,\ot)$, where 
$\Lm$ is a $\bbZ$-lattice, 
to construct the localization functor it is required the ingredients, called 
{\it a real commuting family of graded braiders} $\{(C_i,R_{C_i},\phi_i)\}_{i\in I}$, which consists of 
objects $C_i$ in $\cT$, morphisms $R_{C_i}:(C_i\ot -)\to (-\ot C_i)$ and 
a $\bbZ$-valued function $\phi_i$ on $\Lm$ satisfying
the conditions as in Definition \ref{real-com}. Here note that 
this {\it real commuting family of graded braiders} 
$\{(C_i,R_{C_i},\phi_i)\}_{i\in I}$ plays a similar role as a multiplicative set for the localization of
commutative ring theory.
The resultant localized category $\wtil\cT$ holds several nice properties as in  Theorem \ref{til} and Proposition \ref{til-pro}.

Applying this localization method to the subcategory $\scC_w$ of $R$-gmod associated with
a Weyl group element $w\in W$ one has the localized category $\wtil\scC_w$.  
In particular, since  $\scC_{w_0}$ coincides with the category $R$-gmod 
for finite type $\ge$ and the longest element $w=w_0$  in $W$,  
one can say that through the isomorphism $\cK(\rgmod)\cong 
\cA(\mathfrak n)_{\bbZ[q,q^{-1}]}$ the localized quantum coordinate ring 
${\wtil\aq}_{\bbZ[q,q^{-1}]}$ 
has been defined as $\cK(\tilrgmod)$.

The localized category $\tilrgmod$ possess several interesting properties, such as, 
{\it rigidity} defined as in Definition \ref{rigid} and 
%\begin{df}
%Let $X,Y$ be objects in a monoidal category $\cT$, and $\vep:X\ot Y\to 
%{\bf1}$
%and $\eta:{\bf 1}\to Y\ot X$ morphisms in $\cT$, where $\bf 1$ is 
%the unit object of $\cT$.
%We say that a pair $(X,Y)$ is {\it dual pair} or 
%$X$ is a {\it left dual} to  $Y$ or  
%$Y$ is a {\it right dual} to  $X$ if 
%the following compositions are identities:
%\[
%X\simeq X\ot 1\,\,\mapright{{\rm id}\ot \eta}\,\,X\ot Y\ot X
%\,\,\mapright{\vep\ot {\rm id}}\,\,1\ot X\simeq X,\,\,
%Y\simeq 1\ot Y\,\,\mapright{\eta\ot {\rm id}}\,\,Y\ot X\ot Y
%\,\,\mapright{{\rm id}\ot \vep}\,\,Y\ot 1\simeq Y
%\end{eqnarray*}
%\]
%\end{df}
%\vspace{-20pt}
%We denote a right dual to $X$ by $\cD(X)$ and  a 
%left dual to $X$ by $\cD^{-1}(X)$.
%\vspace{-5pt}
%\begin{thm}[\cite{Loc}]The localized category $\wtil\scC_w$ is left 
%rigid and 
%moreover, for any finite type $\ge$, one has that
% $\tilrgmod$ is { rigid}, i.e.,
%every object in $\tilrgmod$ has left and right duals.
%\end{thm}
the operation $\til\fra$, which is a sort of generalization of the anti-automorphism $\star$ 
on the quantum algebra $\uq=\lan e_i,f_i,q^h\ran_{i\in I,h\in P^*}$ defined by
$e_i^\star=e_i$, $f_i^\star=f_i$ and $(q^h)^\star=q^{-h}$. 
There exists the automorphism $\fra$ on the quiver Hecke algebra $R(\beta)$ as in Sect.9 and then 
it induces the functor $\til\fra:\tilrgmod\to\tilrgmod$ which satisfies 
$\til\fra^2\cong{\rm id}$ and $\til\fra(X\circ Y)\cong \til\fra(Y)\circ \til\fra(X)$ 
for $X,Y\in\tilrgmod$. It follows from the fact $\cK(\tilrgmod)\cong \wtil\aq_{\bbZ[q,q^{-1}]}$ 
that the functor $\til\fra$
defines an anti-automorphism on $\wtil\aq$ and it can be seen as a generalization of 
the  anti-automorphism $\star$.

Furthermore, it is remarkable that for any simple object
$M\in\tilrgmod$ there exist $\lm=\sum_im_i\Lm_i\in P$ and a simple $S\in\rgmod$ such that 
$M\cong C_\lm\circ S$, where $C_\lm\cong C_{1}^{\circ m_1}\circ \cd\circ C_{n}^{\circ m_n}$ up to grading shift.

In the middle of study for the connectedness of the cellular crystal $\bbB_\bfi$
(\cite{Kana-N}), we found the following equality
under some condition in \cite{Kana-N} and without the condition in 
Proposition \ref{pro-bi} below:
\[
\bbB_\bfi=\bigcup_{h\in\cH_\bfi}B^h(\ify)_\bfi,
\]
where $\cH_\bfi$ is a sub $\bbZ$-lattice of 
$\bbB_\bfi$ given as the solution space 
of some system of linear equations \eqref{hi} and 
\[
B^h(\ify)_\bfi:=\{x+h\in \bbZ^N(=\bbB_\bfi)\,|\,x\in B(\ify)_\bfi\}\subset \bbB_\bfi,
\]
which is a subset of $\bbB_\bfi$ obtained by shifting $B(\ify)$ by $h\in\cH_\bfi$.
Note that $\cH_\bfi$ has a natural $\bbZ$-basis $\{\bfh_i\mid i\in I\}$ as in Proposition \ref{basis-h} and 
then any $h\in\cH_\bfi$ can be expressed as $h=\sum_im_i\bfh_i$.
By Lauda-Vazirni results, any simple module in $R$-gmod corresponds to unique 
crystal base element in $B(\ify)$ up to grading shift, where we denote it by 
$\Psi:\{\hbox{simple }S\in\rgmod\}\q\mapright{\sim}\q B(\ify)$.
Here, we observe the following correspondence between the set of 
simple modules in $\tilrgmod$ and 
$\bbB_\bfi$ as
\begin{equation}
\begin{array}{ccc}
\{\hbox{self-dual simple }M\in \tilrgmod\}&\longleftrightarrow &\bbB_\bfi=\bigcup_{h\in\cH_\bfi}B^h(\ify)_\bfi\\
\rotatebox{90}{$\in$}&&\rotatebox{90}{$\in$}\\
M=C_{\sum_im_i\Lm_i}\circ S&\longleftrightarrow& (\sum_im_i\bfh_i)+\Psi(S).
\end{array}
\label{obsv}
\end{equation}
One of the main objectives of this article is to realize this correspondence 
as an isomorphism of crystals. 
In order to pursue it, first we equip the crystal structure with the set 
$\bbB(\tilrgmod)=\{S\mid S\hbox{ is self-dual simple module}\in\tilrgmod\}$
(Theorem \ref{main1}). 
To define the crystal structure on $\bbB(\tilrgmod)$, we shall take into account 
the construction by Lauda-Vazirani\cite{L-V} and consider how we shall 
generalize their construction to the localized category. The most crucial step to 
complete the task is to define the action $\wtil E_i$ on  
$\bbB(\tilrgmod)$
as in \eqref{def-Ei}. It is also essential task to show the well-definedness of
all operators and functions on $\bbB(\tilrgmod)$.
Then, finally, we shall construct the isomorphism of crystals 
$\wtil\Psi:\bbB(\tilrgmod)\mapright{\sim}\bbB_\bfi$ (Theorem \ref{iso}), which seems to be  
a localized version of Lauda-Vazirani's isomorphism from the set of 
simple objects of $\rgmod$ to $B(\ify)\subset \bbB_\bfi$ as mentioned above. 

As has been seen, there exists the functor $\til\fra:\tilrgmod\to\tilrgmod$. 
We shall show the stability of $\bbB(\tilrgmod)$ by $\til\fra$, that is, 
\[
\til\fra(\bbB(\tilrgmod))=\bbB(\tilrgmod),
\]
which seems to be a generalization of the stability $B(\ify)^\star=B(\ify)$ (\cite{K3}).

The organization of the article is as follows.
In Sect.2, we fix the setting of the article and prepare several notation and 
definitions which will be needed later parts.
In Sect.3, we shall review basic definition and properties of crystal bases, crystals and the cellular crystal $\bbB_\bfi$ associated with reduced word of 
Weyl group element $\bfi=i_1\cd i_l$. 
Several crucial properties of cellular crystals will be shown.  
In Sect.4, we shall review basic definitions 
and  several results on quiver Hecke algebra and its modules.
We also review the results by Lauda-Vazirani on categorification of 
crystal $B(\ify)$. Sect.5 is devoted to review 
localization method on monoidal categories following \cite{Loc}.
In Sect.6, the localization method will be applied to the category $\rgmod$
using the {\it determinantial modules} $\{{\bf M}(w_0\Lm_i,\Lm_i)\mid i\in I\}$.
Sect.7 introduces one of the main results of this paper that
the localized monoidal category $\tilrgmod$ possesses the crystal 
$\bbB(\tilrgmod)$. In Sect.8, we shall mention another main results 
of this article that the crystal $\bbB(\tilrgmod)$ is isomorphic to the 
cellular crystal associated with the longest Weyl group element.
In the last section, 
we shall see several miscellaneous results as applications
of the crystal structure on $\bbB(\tilrgmod)$.
It will be shown the stability $\til\fra(\bbB(\tilrgmod))=\bbB(\tilrgmod)$,
which is a generalization of the stability $B(\ify)^\star=B(\ify)$ \cite{K3}.
We also define a certain additive structure on $\bbB(\tilrgmod)$ induced by the one 
of $\bbB_\bfi$ as a free $\bbZ$-lattice.
Finally, several problems will be raised, 
which ask how to translate categorical properties of $\bbB(\tilrgmod)$ to combinatorial ones
of $\bbB_\bfi$, and vice versa.
 Some of them would be discussed elsewhere soon.

\section*{Acknowledgment}
The author would like to thank Masaki Kashiwara from the bottom of his heart
for valuable discussions and helpful advices.
He also thanks Yoshiyuki Kimura and Hironori Oya for their suggestions.

%%%%%%%%%%%%%%%%%%%%%%%%%%%%%%%%%%%
 
\section{Preliminaries}
Let $\ge=\frn\oplus \tt\oplus \frn_-=\lan e_i,h_i,f_i\ran_{i\in I:=\{1,2,\cd,n\}}$
be a simple Lie algebra associated with a Cartan matrix
$ A=(a_{ij})_{i,j\in I}$ where 
$\{e_i,\,f_i,\,h_i\}_{i\in I}$ are the standard Chevalley generators and 
$\frn=\lan e_i\ran_{i\in I}$ (resp. $\tt=\lan h_i\ran_{i\in I}$, 
$\frn_-=\lan f_i\ran_{i\in I}$) is the positive nilpotent subalgebra 
(resp. the Cartan subalgebra, the negative nilpotent subalgebra).

%$\bullet$ $\tt$ : fin.dim.$\bbQ$-vector space and 
%$\tt^*$: dual  vector space of $\tt$.

Let $\{\al_i\}_{i\in I}$ be the set of simple roots of $\ge$ and 
$\lan \q,\q\ran$ a pairing on $\tt \times \tt^*$
satisfying $a_{ij}=(\lan h_i,\al_j\ran)_{i,j\in I}$. We also define 
a symmetric bilinear form $(\,,\,)$ on $\tt^*$ such that
$(\al_i,\al_i)\in 2\bbZ_{> 0}$ and
$\lan h_i,\lm\ran={{2(\al_i,\lm)}\over{(\al_i,\al_i)}}$
for $\lm\in\tt^*$.

Let $P:=\{\lm\in \tt^*\,|\,\lan h_i,\lm\ran\in
\bbZ\,\,\hbox{for any }i\in I\}$ be the weight lattice and 
$P_+:=\{\lm\in P\mid \lan h_i,\lm\ran\geq0\hbox{ for any }i\in I\}$ the set of dominant weights.
Set $Q:=\oplus_{i\in I}\bbZ\al_i$ (resp. $Q_+:=\sum_{i\in I}\bbZ_{\geq0}\al_i$), which is called the root lattice (resp. positive root lattice). 
For an element  $\beta=\sum_i m_i\al_i\in Q_+$ define $|\beta|=\sum_i m_i$, which is called the height of $\beta$.
Let $W=\lan s_i\,|\,s_i\ran_{i\in I}$ be the Weyl group associated with $P$,
where $s_i$ is the simple reflection defined by $s_i(\lm)
=\lm-\lan h_i,\lm\ran\al_i$ $(\lm\in P)$.

We denote the dual weight lattice of $P$ by
$P^*:=\{h\in \tt\,|\,\lan h,P\ran\subset\bbZ\}$. 
Let $ \uq:=\lan e_i,\,f_i, \,q^h \ran_{i\in I,h\in P^*}$ be the 
quantum algebra associated with $\ge$ with the defining relations
(see e.g.,\cite{K1,K3}) and $\uqm:=\lan f_i\ran_{i\in I}$ (resp. 
$\uqp:=\lan e_i\ran_{i\in I}$) the 
negative (resp. positive) nilpotent subalgebras of $\uq$.
We also define the $\bbZ$-form $U^-_{\bbZ[q,q^{-1}]}(\ge)$ of $\uqm$ as in 
\cite{Loc}.   

Now, let us define the 
{\it (unipotent) quantum coordinate ring $\aq$} by
\[
\cA_q(\frn)=\bigoplus_{\beta\in Q_-}\cA_q(\frn)_\beta\qq
\cA_q(\frn)_\beta:={\rm Hom}_{\bbQ(q)}(U^+_q(\ge)_{-\beta},\bbQ(q))
\]
Note that $\uqm\cong \aq$ as a $\bbQ(q)$-algebra. The $\bbZ$-form 
$\cA(\mathfrak n)_{\bbZ[q,q^{-1}]}$ is defined as in \cite{Loc}.
 
%%%%%%%%%%%%%%%%%%%%%%%%%%%%%%%%%%%
%%%%%%%%%%%%%%%%%%%%%
\section{Crystal Bases and Crystals}
%%%%%%%%%%%%%%%%%
\subsection{Crystal Base of $\uqm\cong \aq$}
%$\bl$ $V(\lm)$ be the irr. highest weight mod. with h.w.v $u_\lm$ 
%($u_\lm\in P_+)$. \\
%\medskip
As in \cite{K1}, there exists the crystal base $(L(\ify),B(\ify))$ of $\uqm$ defined by
\begin{eqnarray*}
&& L(\ify):=\sum_{k\geq0,i_1,\cd,i_k\in I}\bbA\til f_{i_1}\cd\til
 f_{i_k}u_\ify,\\
&&B(\ify)=\{\til f_{i_1}\cd\til
 f_{i_k}u_\ify\,{\rm mod}\, q L(\ify)\,|\,k\geq0, i_1,\cd,i_k\in I\}\setminus\{0\},\\
%&&\\
&&\vep_i(b)={\rm max}\{k:\eit^kb\ne0\},\q
\vp_i(b)=\vep_i(b)+\lan h_i,\wt(b)\ran,
\end{eqnarray*}
where $u_\ify=1\in \uq$, 
$\eit$ and $\fit\in{\rm End}_{\bbQ(q)}(\uq)$ are the {\it Kashiwara operators}
(\cite{K1}) and 
$\bbA\subset \bbQ(q)$ is the local subring at $q=0$.

%%%%%%%%%%%%%%%%%%%%%%%%%
\subsection{Crystals}
%Let ${}^t\bbZ:=\bbZ\sqcup\{-\ify\}$ be a tropical semi-field with the 
%multiplication $+$ and the  summation $\max/\min$.
We define the notion of {\it crystal} as in \cite{K3}, which is 
a combinatorial object abstracting the properties of crystal bases:
\begin{df}[\cite{K3}]\label{cryst}
A 6-tuple $(B,\wt, \{\vep_i\},\{\vp_i\}, \{\eit\},\{\fit\})_{i\in I}$ is 
a {\it crystal} if $B$ is a set and there exists a certain special element 
$0$ outside of $B$ and maps:
\begin{eqnarray}
&&{\rm wt}:B\to P,\label{wtp-c}\q
\vep_i:B\to\bbZ\sqcup\{-\ify\},\q\vp_i:B\to\bbZ\sqcup\{-\ify\}\q\,(i\in I),
\\
&&\eit:B\sqcup\{0\}\to B\sqcup\{0\},\q
\fit:B\sqcup\{0\}\to B\sqcup\{0\}\,\,(i\in I),\label{eitfit-c}
\end{eqnarray}
satisfying :
\begin{enumerate}
\item
$\vp_i(b)=\vep_i(b)+\lan h_i,\wt(b)\ran$.
\item
If  $b,\eit b\in B$, then $\wt(\eit b)=\wt(b)+\al_i$, 
$\vep_i(\eit b)=\vep_i(b)-1$, $\vp_i(\eit b)=\vp_i(b)+1$.
\item
If $b,\fit b\in B$, then $\wt(\fit b)=\wt(b)-\al_i$, 
$\vep_i(\fit b)=\vep_i(b)+1$, $\vp_i(\fit b)=\vp_i(b)-1$.
\item
For $b,b'\in B$ and $i\in I$, one has 
$\fit b=b'$ iff $b=\eit b'.$
\item
If $\vp_i(b)=-\ify$ for $b\in B$, then $\eit b=\fit b=0$
and $\eit(0)=\fit(0)=0$.
\end{enumerate}
Here, a {\it crystal graph} of crystal $B$ is a $I$-colored 
oriented graph defined by $b\mapright{i} b'\Leftrightarrow \fit(b)=b'$ for $b,b'\in B$.
\end{df}

\begin{df}[\cite{K3}]\label{morph}
For crystals $B_1$ and $B_2$, 
$\Psi$ is a {\it strict embedding} (resp. {\it isomorphism})
from $B_1$ to $B_2$ if 
$\Psi:B_1\sqcup\{0\}\to B_2\sqcup\{0\}$ is an {injective }
(resp. {bijective}) map satisfying that $\Psi(0)=0$, 
$\wt(\Psi(b))=\wt(b)$, 
$\vep_i(\Psi(b))=\vep_i(b)$ and 
$\vp_i(\Psi(b))=\vp_i(b)$ for any $b\in B_1$ and 
$\Psi$ commutes with 
all $\eit$'s and $\fit$'s,.
\end{df}

We obtain the tensor structure of crystals as follows(\cite{K1,K3}):
\begin{pro}\label{tensor}
For crystals $B_1$ and $B_2$,  
set 
\[
B_1\ot B_2=\{b_1\ot b_2:=(b_1,b_2)\mid b_1\in B_1 ,\, b_2\in B_2\}
(=B_1\times B_2).
\]
Then, $B_1\ot B_2$ becomes a crystal by defining:
\begin{eqnarray}
&&wt(b_1\ot b_2)=wt(b_1)+wt(b_2),\\
&&\vep_i(b_1\ot b_2)={\hbox{max}}(\vep_i(b_1),
  \vep_i(b_2)-\lan h_i,wt(b_1)\ran),
\label{tensor-vep}\\
&&\vp_i(b_1\ot b_2)={\hbox{max}}(\vp_i(b_2),
  \vp_i(b_1)+\lan h_i,wt(b_2)\ran),
\label{tensor-vp}\\
&&\eit(b_1\ot b_2)=
\left\{
\begin{array}{ll}
\eit b_1\ot b_2 & {\mbox{ if }}\vp_i(b_1)\geq \vep_i(b_2)\\
b_1\ot\eit b_2  & {\mbox{ if }}\vp_i(b_1)< \vep_i(b_2),
\end{array}
\right.
\label{tensor-e}
\\
&&\fit(b_1\ot b_2)=
\left\{
\begin{array}{ll}
\fit b_1\ot b_2 & {\mbox{ if }}\vp_i(b_1)>\vep_i(b_2)\\
b_1\ot\fit b_2  & {\mbox{ if }}\vp_i(b_1)\leq \vep_i(b_2).
\label{tensor-f}
\end{array}
\right.
\end{eqnarray}
\end{pro}

%%%%%%%%%%%%%%%%%%%%%
% 
%\section{Strict Embedding of Crystals}
%\begin{df}
%A crystal $B$ is a {upper(resp. lower) semi-normal crystal} if 
%for any $i\in I$, we have 
%$\vep_i(b)=\max\{k\,|\,\eit^k(b)\ne0\}$ (resp. 
%$\vp_i(b)=\max\{k\,|\,\fit^k(b)\ne0\}$). 
%If a crystal $B$ is upper and lower semi-noraml, it is called 
%a { normal crystal}.
%\end{df}
% {}
% In this sense, the crystal $B(\lm)$ $(\lm\in P_+)$ is a 
%normal crystal and $B(\ify)$ is a upper semi-normal crystal.
% 
%
%
%
%\begin{df}\label{morph}
%For crystals $B_1,\, B_2$, 
%$\Psi$ is a {\bf strict embedding} (resp. {\bf isomorphism})
%from $B_1$ to $B_2$ if 
%$\Psi:B_1\sqcup\{0\}\to B_2\sqcup\{0\}$ is a {injective }
%(resp. {bijective})map satisfying : 
%\begin{enumerate}
%\item
%$\Psi(0)=0$.
%\item
%For $b\in B_1$, $\wt(\Psi(b))=\wt(b)$, 
%$\vep_i(\Psi(b))=\vep_i(b)$ and 
%$\vp_i(\Psi(b))=\vp_i(b)$.
%\item
%$\Psi:B_1\to B_2$ commutes with 
%all $\eit$'s and $\fit$'s,.
%\end{enumerate}
%\end{df}
%%Let us define the category of crystals by taking morphisms of crystals 
%%as morphisms. It is  denoted by $\cC\cR$.
%
% 
%%%%%%%%%%%%%%%%%%%%%%

\begin{ex}\label{ex-cry}
For $i\in I$, set
$ B_i:=\{(n)_i\,|\,n\in \bbZ\}$ and 
\begin{eqnarray*}
&&\wt((n)_i)=n\al_i,\,\,
\vep_i((n)_i)=-n,\,\,
\vp_i((n)_i)=n,\,\,\label{bidata}\\
&&
\vep_j((n)_i)=\vp_j((n)_i)=-\ify\,\,\,(i\ne j),\label{biji}\\
&&
\eit((n)_i)=(n+1)_i,\q\fit((n)_i)=(n-1)_i,\label{bieitfit}\\
&&
\til e_j((n)_i)=\til f_j((n)_i)=0\,\,\,(i\ne j).
\label{biejfj}
\end{eqnarray*}
\end{ex}
Then $B_i$ ($i\in I$) possesses a crystal structure.
Note that as a set the crystal $B_i$ can be identified with the set of integers 
$\bbZ$.
\subsection{Explicit structure of the crystal
$B_{i_1}\ot\cd\ot B_{i_m}$}
In this subsection, we shall investigate explicit structure of 
a tensor product of the crystal $B_i$'s.

Fix a sequence of indices $\bfi=(i_1,\cd,i_m)\in I^m$
and write 
\[
(x_1,\cd,x_m):=\til f_{i_1}^{x_1}(0)_{i_1}
\ot\cd\ot \til f_{i_m}^{x_m}(0)_{i_m}=(-x_1)_{i_1}\ot\cd\ot(-x_m)_{i_m},
\]
where if $n<0$, then $\fit^n(0)_i$ means $\eit^{-n}(0)_i$.
Note that here we do not necessarily assume that $\bfi$ is a reduced word
though later we will take $\bfi$ to be a reduced longest word.
%Set
%\begin{equation}
%\ZZ^{\ify}
%:=\{(\cd,x_k,\cd,x_2,x_1): x_k\in\ZZ
%\,\,{\rm and}\,\,x_k=0\,\,{\rm for}\,\,k\gg 0\};
%\label{uni-cone}
%\end{equation}
By the tensor structure of crystals 
in Proposition \ref{tensor}, for the sequence  $\bfi$ as above,  we can describe 
the explicit crystal structure on $\bbB_\bfi:=B_{i_1}\ot \cd\ot B_{i_m}$ 
as follows: 
For $x=(x_1,\cd,x_m)\in \bbB_\bfi$, define
\begin{equation*}
 \sigma_k(x):=x_k+\sum_{j<k}\lan h_{i_k},\al_{i_j}\ran x_j
\label{sigma-def}
\end{equation*}
and for $i\in I$ define
\begin{eqnarray*}
&&\TY(\wtil\sigma,i,)(x):=\max\{\sigma_k(x)\,|\,1\leq k\leq m\,{\rm and}\,
i_k=i\},\label{simgak}\\
&&\TY(\wtil M,i,)=\TY(\wtil M,i,)(x):=\{k\,|\,1\leq k\leq m,\,i_k=i,\,
\sigma_k(x)=\TY(\wtil\sigma,i,)(x)\},
\label{Mi}
\\
&&
\TY(\wtil m,i,f)=\TY(\wtil m,i,f)(x):=\max\,\TY(\wtil M,i,)(x),\q
\TY(\wtil m,i,e)=\TY(\wtil m,i,e)(x):=\min\,\TY(\wtil M,i,)(x).
\label{mfme}
\end{eqnarray*}
Now, the actions of the Kashiwara operators $\eit,\fit$ and the functions
$\vep_i,\vp_i$ 
and $\wt$  are written explicitly:
\begin{eqnarray}
&& \fit(x)_k:=x_k+\del_{k,\TY(\wtil m,i,f)},\qq\qq
%\TY(\wtil\sigma,i,)(x)>0\,{\rm then}\,
\eit(x)_k:=x_k-\del_{k,\TY(\wtil m,i,e)},\\
&&
\wt(x):=-\sum_{k=1}^m x_k\al_{i_k},\q
\vep_i(x):=\TY(\wtil\sigma,i,)(x),\q
\vp_i(x):=\lan h_i,\wt(x)\ran+\vep_i(x).
\label{wt-vep-vp-ify}
\end{eqnarray}

Define the function $\beta^{(\bfi)}_k$ on $\bbB_\bfi$ by :
\begin{equation}
\beta^{(\bfi)}_k(x):=\sigma_{k^+}(x)-\sigma_{k}(x)
=x_k+\sum_{k<j<k^+}\lan h_i,\al_{i_j}\ran x_j+x_{k^+},
\label{betak}
\end{equation}
for $x=(x_1,\cd,x_m)\in \bbB_\bfi$, where 
for $k\in[1,m]$, $k^+$ (resp. $k^-$) is the minimum (resp. maximum) number $j\in[1,m]$ such that
$k<j$ (resp. $l<k$) and $i_k=i_j$ if it exists, otherwise $m+1$ (resp. 0).
Here one knows that $\TY(\wtil m,i,f)(x)$ and $\TY(\wtil m,i,e)(x)$ are determined 
by $\{\beta^{(\bfi)}_k(x)\,|1\leq  k\leq m,\,\,i_k=i\}$.
%%%%%%%%%%%%%%%%%%%%%%%%%%%%%%%%%%%
\subsection{Braid-type isomorphism}
Set 
$c_{ij}:=\lan h_i,\al_j\ran\lan h_j,\al_i\ran$,
$c_1:=-\lan h_i,\al_j\ran$ and $c_2:=-\lan h_j,\al_i\ran$.
In the sequel, for $x\in \bbZ$, put
\[x_+:=\begin{cases}
x& \hbox{ if }x\geq 0,\\
0& \hbox{ if }x<0.
\end{cases}
\]

\begin{pro}[\cite{N1}]
There exist the following 
isomorphisms of crystals $\phi^{(k)}_{ij}$ ($k=0,1,2,3$)
\begin{enumerate}
\item
If $c_{ij}=0$,
\beq
\phi^{(0)}_{ij}:B_i\ot B_i\mapright{\sim} B_j\ot B_i,\label{bt0}
\eeq
where $\phi^{(0)}_{ij}((x)_i\ot(y)_j)=(y)_j\ot (x)_i$.
\item
If $c_{ij}=1$,
\beq
\phi^{(1)}_{ij}:B_i\ot B_j\ot B_i\mapright{\sim} B_j\ot B_i\ot B_j,
\label{bt1}
\eeq
where 
$$
\phi^{(1)}_{ij}((x)_i\ot(y)_j\ot(z)_i)=
(z+(-x+y-z)_+)_j\ot (x+z)_i\ot (y-z-(-x+y-z)_+)_j.
$$
\item
If $c_{ij}=2$,
\beq
\phi^{(2)}_{ij}:B_i\ot B_j\ot B_i\ot B_j
\mapright{\sim} B_j\ot B_i\ot B_j\ot B_i,
\label{bt2}
\eeq
where 
$\phi^{(2)}_{ij}$ is given by the following:
for $(x)_i\ot(y)_j\ot(z)_i\ot (w)_j$ we set 
%%% $P:=-a+c_1b-c$, $Q:=-b+c_2c-d$ and $R:=-a+c-c_1d$, and 
$(X)_j\ot (Y)_i\ot (Z)_j\ot (W)_i
:=\phi^{(2)}_{ij}((x)_i\ot(y)_j\ot(z)_i\ot (w)_j)$.
\beqn
X & = & w+(-c_2x+y-w+c_2(x-c_1y+z)_+)_+,
\label{X}\\
Y & = & x+c_1w+(-x+z-c_1w+(x-c_1y+z)_+)_+,
\label{Y}\\
Z & = & y-(-c_2x+y-w+c_2(x-c_1y+z)_+)_+,
\label{Z}\\
W & = & z-c_1w-(-x+z-c_1w+(x-c_1y+z)_+)_+.
\label{W}
\eeqn
\item
If $c_{ij}=3$, the map 
\beq
\phi^{(3)}_{ij}:B_i\ot B_j\ot B_i\ot B_j\ot B_i\ot B_j
\mapright{\sim} B_j\ot B_i\ot B_j\ot B_i\ot B_j\ot B_i,
\label{bt3}
\eeq
is defined by the following: for 
$(x)_i\ot (y)_j\ot (z)_i\ot (u)_j\ot (v)_i\ot (w)_j$
we set 
$A:=-x+c_1y-z$,
$B:=-y+c_2z-u$,
$C:=-z+c_1u-v$ and 
$D:=-u+c_2v-w$.
Then $(X)_j\ot (Y)_i\ot (Z)_j\ot (U)_i\ot (V)_j\ot (W)_i
:=\phi^{(3)}_{ij}((x)_i\ot (y)_j\ot (z)_i\ot (u)_j\ot (v)_i\ot (w)_j)$
is given by
\begin{eqnarray*}
X & = & w+(D+(c_2C+(2B+A_+)_+)_+)_+,\\
Y & = & x+c_1w+(c_1D+(3C+(2c_1B+2A_+)_+)_+)_+,\\
Z & = & y+u+w-X-V,\\
U & = & x+z+v-Y-W,\\
V & = & u-w-(2D+(2c_2C+(3B+c_2A_+)_+)_+)_+,\\
W & = & v-c_1w-(c_1D+(2C+(c_1B+A_+)_+)_+)_+.
\end{eqnarray*}
\end{enumerate}
They also satisfy $\phi_{ij}^{(k)}\circ \phi_{ji}^{(k)}={\rm id}$.

We call such isomorphisms of crystals  {\it braid-type isomorphisms}.
\end{pro}

We also define a {\it braid-move} on the set of reduced words of $w\in W$
to be a composition of the following transformations 
induced from the braid relations:
\begin{eqnarray*}
&&\cd ij\cd\to \cd ji\cd (c_{ij}=0),\q
\cd iji\cd\to \cd jij\cd (c_{ij}=1),\\
&&\cd ijij\cd\to \cd jiji\cd (c_{ij}=2),\q
\cd ijijij\cd\to \cd jijiji\cd (c_{ij}=3),
\end{eqnarray*}
which are called 2-move, 3-move, 4-move, 6-move respectively.

%%%%%%%%%%%%%%%%%%%%%%%%%%%%%%%%%%%%%%%%%%%%%%%%%
\subsection{Cellular Crystal $\bbB_\bfi=\bbB_{i_1i_2\cd i_k}
=B_{i_1}\ot\cd\ot B_{i_k}$}
For a reduced word $\bfi=i_1i_2\cd i_k$ of some Weyl group element, 
we call the crystal $\bbB_\bfi:=B_{i_1}\ot\cd\ot B_{i_k}$ 
a {\it cellular crystal}
associated with a reduced word $\bfi$. Indeed, 
it is obtained by applying the tropicalization functor to the 
geometric crystal on the Langlands-dual { Schubert cell ${}^LX_w$}, where 
$w=s_{i_1}\cd s_{i_k}$ is an element of the Well group $W$ (\cite{N3}).
%We obtain the following in [KN]:
It is immediate from the braid-type isomorphisms that 
for any $w\in W$ and its arbitrary reduced words $i_1\cd i_l$ and 
$j_1\cd j_l$, we get the following isomorphism of crystals:
\begin{equation}
B_{i_1}\ot\cd\ot B_{i_l}\cong B_{j_1}\ot\cd\ot B_{j_l}.
\label{iso-cellular}
\end{equation}

%%%%%%%%%%%%%%%%%%%%%%%%%%%%%%%%%%%%%%%%%%%
\subsection{Half potential and the crystal $B(\ify)$}\label{subsec-half}
For a Laurent polynomial $\phi(x_1,\cd,x_n)$ with positive coefficients, 
the tropicalization of $\phi$ is denoted by $\wtil\phi
:={\rm Trop}(\phi)$, which is given by the rule:
${\rm Trop}(ax+by)={\rm min}(x,y)$ with $a,b>0$, 
${\rm Trop}(xy)=x+y$ and ${\rm Trop}(x/y)=x-y$ and ${\rm Trop}(c)=0$ for $c>0$.
In \cite{Kana-N}, the crystal $B(\ify)$ has been realized as a certain subset of 
$\bbB_\bfi$ defined as follows:
\begin{thm}[{\cite[Theorem 5.11]{Kana-N}}]
Define the subset of $\bbB_\bfi$:
\[
(\wtil \bbB^-_{w_0})_{\Phi^{(+)},\Theta_\bfi}
=\{x=(x_1,\cd,x_N)\in \bbB_\bfi\mid
\wtil \Phi^{(+)}(x)\geq0\},
\]
where $\bbB^-_{w_0}:=B^-\cap B\ovl{w_0}B$ holds a certain geometric crystal structure,  
$\wtil\Phi^{(+)}$ is a tropicalization of the half potential
$\Phi^{(+)}$ which is a  
Laurent polynomial with positive coefficients in $N$ variables and 
$\Theta_\bfi$ is a certain positive structure on the geometric crystal 
$\bbB^-_{w_0}$.
Then, $(\wtil \bbB^-_{w_0})_{\Phi^{(+)},\Theta_\bfi}\cong B(\ify)$.
\end{thm}
\begin{rem}\label{trop-real}
To define the crystal structure on 
$(\wtil \bbB^-_{w_0})_{\Phi^{(+)},\Theta_\bfi}$,
it is supposed that if $\eit x\not\in (\wtil \bbB^-_{w_0})_{\Phi^{(+)},\Theta_\bfi}$, 
then $\eit x=0$. Thus, in this sense, the embedding 
$B(\ify)\cong 
(\wtil \bbB^-_{w_0})_{\Phi^{(+)},\Theta_\bfi}\hookrightarrow \bbB_\bfi$ is not
a strict embedding. In \cite{N3,N-Z}, it has been given the strict embedding of 
$B(\ify)\hookrightarrow \bbB_\bfi$, which is called "Kashiwara embedding" and 
the method to describe the image of this embedding is called "polyhedral realization".
\end{rem}

%%%%%%%%%%%%%%%%%%%%%%%%%%%%%%%%%%%%%%%%%%%%
\subsection{Subspace $\cH_\bfi$}

Fix a reduced longest word $\bfi=i_1\cd i_N$ and 
take the function 
$\beta^{(\bfi)}_k(x)=x_k+\sum_{k<j<k^+}\lan h_{i_k}\al_{i_j}\ran x_j+x_{k^+}
\q (1\leq  k\leq N)$ as in \eqref{betak}.
In what follows, let us identify the $\bbZ$-lattice $\bbZ^N$ with $B_\bfi$ and then 
we define the summation of elements $x=(x_1,\cd,x_N)$ and $y=(y_1,\cd,y_N)$
by $x+y=(x_1+y_1,\cd,x_N+y_N)$ as a standard one in $\bbZ^N$.
Here, we define the subspace $\cH_\bfi\subset\bbZ^N$ by 
\begin{equation}
\cH_\bfi:=\{x\in \bbZ^N(=\bbB_\bfi)\,|\,\beta^{(\bfi)}_k(x)=0
\hbox{ for any } k{\hbox{ such that }}k^+\leq N\}\subset \bbB_\bfi.
\label{hi}
\end{equation}
 
%For $\bfi$, set $\color{navyblue}\TY(\al,k,)
%=s_{i_N}s_{i_{N-1}}\cd s_{i_{k+1}}(\al_{i_k})$ 
%and $\color{navyblue}\Del:=\{\TY(\al,k,)\,|\,1\leq k\leq N\}$. 
%It is well-known that
%$\Del$ is the set of all positive roots. 
The following proposition describes the  result of \cite[Lemma 9.1]{Kana-N} explicitly:
\begin{pro}[\cite{Kana-N}]\label{basis-h}
For $\bfi=i_1i_2\cd i_N$, $k=1,2,\cd, N$ 
 and a fundamental weight $\Lm_i$, set 
\begin{equation}
\TY(h,k,i):=\lan h_{i_k},s_{i_{k+1}}\cd s_{i_N}\Lm_i\ran\q\hbox{ and }\q
\bfh_i:=(\TY(h,1,i),\TY(h,2,i),\ldots,\TY(h,N,i))\in \bbB_\bfi
\label{hik}
\end{equation}
Then, we obtain that $\{\bfh_1,\cd,\bfh_n\}$ is a $\bbZ$-basis of $\cH_\bfi$, 
namely,
\begin{equation}
\cH_\bfi=
\bbZ\bfh_1\oplus\bbZ\bfh_2\oplus\cd\oplus \bbZ\bfh_n.
%\left\{\left(\sum_i\TY(h,1,i)h_i,\sum_i\TY(h,2,i)h_i,\cd,
%\sum_i\TY(h,N,i)h_i\right)\,|\,(h_1,h_2,\cd,h_n)\in\bbZ^n\right\}
%\label{Hm}
\end{equation}
\end{pro}
{\sl Proof.}
Let $\{\al'_i\}_i$, $\{h'_i\}_i$ and $\{s'_i\}_i$ be the simple roots, the simple co-roots
and the simple reflections of the Langlands dual Lie algebra $\ge^\vee$ respectively. Define 
$m^{(k)}_i\in\bbZ_{\geq 0}$ $(k\in [1,N],\, i\in I)$ by 
\[
{\al'}^{(k)}:=s'_{i_N}s'_{i_N-1}\cd s'_{k+1}(\al'_{i_k})=\sum_{i\in I}m^{(k)}_i\al'_i.
\]
By \cite[Lemma 9.1]{Kana-N}, one has that 
$\{{\bf m}_i:=(m^{(1)}_i,m^{(2)}_i,\cd,m^{(N)}_i)\mid i\in I\}$ is a $\bbZ$-basis of $\cH_\bfi$. 
Thus, it suffices to show that 
$h^{(k)}_i=m^{(k)}_i$ for any $k\in[1,N]$ and $i\in I$.

Let us define the set of paths from $a$ to $b$ $(a,b\in \bbZ, \,a\geq b)$ by
\[
\cP(a,b):=\{(a,j_1,j_2,\ldots,j_l,b)\mid a>j_1>j_2>\cdots >j_l>b,\,l\geq 0\},
\]
where set $\cP(a,a)=\emptyset$ and $l=-1$.
The following lemma is obtained by direct calculations.
\begin{lem}\label{explicit-h-m}
We obtain the following explicit formulas:
\begin{eqnarray}
&&\lan h_{i_k},s_{i_{k+1}}\cd s_{i_{p-1}}(\al_{i_p})\ran\label{h-al}\\
&&\qq=\sum_{(p,j_1,\cd,j_l,k)\in \cP(p,k)}
(-1)^l\lan h_{i_k},\al_{i_{j_l}}\ran\lan h_{i_{j_l}},\al_{i_{j_{l-1}}}\ran
\cd\lan h_{i_{j_2}},\al_{i_{j_1}}\ran\lan h_{i_{j_1}},\al_{i_p}\ran\q(p>k), \nn \\
&&s'_{i_N}s'_{i_{N-1}}\cd s'_{i_{k+1}}(\al'_{i_k}) \label{m-al}\\
&&\q =\sum_{L=k}^N\sum_{\,\,(L,j_1,\cd,j_l,k)\in\cP(L,k)}
(-1)^{l+1}\lan h'_{i_L},\al'_{i_{j_1}}\ran\lan h'_{i_{j_1}},\al'_{i_{j_2}}\ran\cd
\lan h'_{i_{j_{l-1}}},\al'_{i_{j_l}}\ran\lan h'_{i_{j_l}},\al'_{i_k}\ran\al'_{i_L},\nn
\end{eqnarray}
where note that in \eqref{m-al} if $k=L$, namely $\cP(L,k)=\emptyset$, then the corresponding term is $\al'_{i_k}$.
\end{lem}
For $i\in I$ and $i_k,\ldots,i_N$, let $J:=\{j,j^+,j^{2+},\cd,j^{r+}\}$ be the set of all indices 
such that $i=i_j=\cd=i_{j^{r+}}$ and $k<j<j^+<\cd< j^{r+}\leq N$. If $J=\emptyset$,  then
we get $\lan h_{i_k},s_{i_{k+1}}\cd s_{i_N}(\Lm_i)\ran=\del_{i,i_k}$.
Assume that $J\ne \emptyset$. 
Then, by \eqref{h-al} one gets
\begin{eqnarray*}
&&\lan h_{i_k},s_{i_{k+1}}\cd s_{i_N}(\Lm_i)\ran=
\lan h_{i_k},s_{i_{k+1}}\cd s_{i_{j^{r+}}}(\Lm_i)\ran
=\lan h_{i_k},s_{i_{k+1}}\cd s_{i_{j^{r+}-1}}(\Lm_i-\al_i)\ran\\
&&\q=\lan h_{i_k},s_{i_{k+1}}\cd s_{i_{j^{(r-1)+}-1}}(\Lm_i)\ran 
-\lan h_{i_k},s_{i_{k+1}}\cd s_{i_{j^{(r-1)+}-1}}(\al_i)\ran
-\lan h_{i_k},s_{i_{k+1}}\cd s_{i_{j^{r+}-1}}(\al_i)\ran=\cd\\
&&\q=\del_{i,i_k}-\sum_{p=1}^r
\lan h_{i_k},s_{i_{k+1}}\cd s_{i_{j^{p+}-1}}(\al_i)\ran\\
&&\q=\del_{i,i_k}+\sum_{p=1}^r
\sum_{(j^{p+},j_1,\cd,j_l,k)\in \cP(j^{p+},k)}
(-1)^{l+1}\lan h_{i_k},\al_{i_{j_l}}\ran\lan h_{i_{j_l}},\al_{i_{j_{l-1}}}\ran
\cd\lan h_{i_{j_2}},\al_{i_{j_1}}\ran\lan h_{i_{j_1}},\al_{i}\ran
\end{eqnarray*}
On the other-hand, by \eqref{m-al} one gets
\begin{equation}
\hspace{20pt} m^{(k)}_i=\del_{i,i_k}+\sum_{p=1}^r
\sum_{\,\,(j^{p+},j_1,\cd,j_l,k)\in\cP(j^{p+},k)}
(-1)^{l+1}\lan h'_i,\al'_{i_{j_1}}\ran\lan h'_{i_{j_1}},\al'_{i_{j_2}}\ran\cd
\lan h'_{i_{j_{l-1}}},\al'_{i_{j_l}}\ran\lan h'_{i_{j_l}},\al'_{i_k}\ran,
\label{mik} 
\end{equation}
where note that 
the term $\del_{i,i_k}$ in \eqref{mik} is derived from $\cP(k,k)$ in \eqref{m-al}.
Comparing these formulae and using the fact $\lan h_i,\al_j\ran=\lan h'_j,\al'_i\ran$, we obtain $m^{(k)}_i= 
\lan h_{i_k},s_{i_{k+1}}\cd s_{i_N}(\Lm_i)\ran=h^{(k)}_i$.
\qed
%%%%%%%%%%%%%%%%%
 
\begin{ex}\label{G2}
In $\ge=G_2$-case. Set $a_{12}=-1$ and $a_{21}=-3$. 
Taking a reduced longest word $\bfi=121212$, one has
\[
\beta^{(\bfi)}_1(x)=x_1-x_2+x_3,\q
\beta^{(\bfi)}_2(x)=x_2-3x_3+x_4,\q
\beta^{(\bfi)}_3(x)=x_3-x_4+x_5,\q
\beta^{(\bfi)}_4(x)=x_4-3x_5+x_6.
%\beta_k(x)=\begin{cases}x_k-x_{k+1}+x_{k+2}&k=1,3\\
%x_k-3x_{k+1}+x_{k+2}&k=2,4
%stuple\end{cases}
\]
By the formula \eqref{hik}, one gets
\[
\bfh_1=(1,3,2,3,1,0),\qq
\bfh_2=(0,1,1,2,1,1).
\]
%\begin{eqnarray*}
%&&\TY(\al,1,)={\al_2}, \q
%\TY(\al,2,)=s_2(\al_1)={\al_1}+{3\al_2},\q
%\TY(\al,3,)=s_2s_1(\al_2)={\al_1}+{2\al_2},\\
%&&\TY(\al,4,)=s_2s_1s_2(\al_1)={2\al_1}+{3\al_2},\q
%\TY(\al,5,)=s_2s_1s_2s_1(\al_2)={\al_1}+{\al_2},\q
%\TY(\al,6,)={\al_1}
%\end{eqnarray*}
Then the solution space $\cH_\bfi$ 
of $\beta^{(\bfi)}_1(x)=\beta^{(\bfi)}_2(x)=\beta^{(\bfi)}_3(x)=\beta^{(\bfi)}_4(x)=0$
is given  by 
\[
\cH_\bfi=\{{ c_1}\bfh_1+{ c_2}\bfh_2=
({ c_1},{ c_2}+{3c_1}, { c_2}+{2c_1},
{ 2c_2}+{3c_1},
{ c_2}+{ c_1},{ c_2})\mid
c_1,c_2\in\bbZ\}.
\]
\end{ex}
%The explicit forms of $\TY(\al,k,)$'s:
\begin{lem}\label{310}
The braid-type isomorphisms are well-defined on $\cH_\bfi$, that is, 
$\phi_{ij}^{(k)}(\cH_\bfi)=\cH_{\bfi'}$, where $\bfi'$ is the  
reduced word obtained 
by applying the corresponding braid-moves. We also obtain the following 
formula:
\begin{enumerate}
\item
Case $c_{ij}=0$: For any $h=(\cd,x,y,\cd)=\cd\ot(-x)_i\ot (-y)_j\ot\cd
\in \cH_\bfi$, 
applying the braid-type isomorphism $\phi_{ij}^{(0)}$ on $(x,y) $ in $h$,
we have 
\begin{equation}
\phi^{(0)}_{ij}(h)=(\cd,y,x,\cd)=\cd\ot(-y)_j\ot (-x)_i\ot\cd
\in\cH_{\bfi'}
\label{br0}
\end{equation} 
\item
Case $c_{ij}=1$: For any $h=(\cd,x,y,z,\cd)=\cd\ot(-x)_i\ot (-y)_j\ot(-z)_i\ot\cd
\in \cH_\bfi$,  
applying the braid-type isomorphism $\phi_{ij}^{(1)}$ on $(x,y,z) $ in $h$,
we have 
\begin{equation}
\phi^{(1)}_{ij}(h)=(\cd,z,y,x,\cd)=\cd\ot(-z)_j\ot (-y)_i\ot(-x)_j\ot\cd
\in\cH_{\bfi'}
\label{br1}
\end{equation} 
\item
Case $c_{ij}=2$: For $h=(\cd,x,y,z,w,\cd)=\cd\ot(-x)_i\ot (-y)_j\ot(-z)_i\ot(-w)_j\cd
\in \cH_\bfi$, 
applying the braid-type isomorphism $\phi_{ij}^{(2)}$ on $(x,y,z,w)$ in $h$,
we have 
\begin{equation}
\phi^{(2)}_{ij}(h)=(\cd,w,z,y,x,\cd)
=\cd\ot(-w)_j\ot(-z)_i\ot(-y)_j\ot(-x)_i\ot\cd
\in\cH_{\bfi'}
\label{br2}
\end{equation} 
\item
Case $c_{ij}=3$: For $h=(\cd,x,y,z,u,v,w,\cd)=
\cd\ot(-x)_i\ot (-y)_j\ot(-z)_i\ot(-u)_j\ot(-v)_i\ot(-w)_j\cd
\in \cH_\bfi$, 
applying the braid-type isomorphism $\phi_{ij}^{(3)}$ on $(x,y,z,u,v,w)$ in $h$,
we have 
\begin{equation}
\phi^{(3)}_{ij}(h)=(\cd,w,v,u,z,y,x,\cd)
=\cd\ot(-w)_j\ot (-z)_i\ot(-y)_j\ot(-x)_i\ot\cd
\in\cH_{\bfi'}
\label{br3}%see mature-note pp48
\end{equation} 
\end{enumerate}
\end{lem}
{\sl Proof.}
The formula \eqref{br0}--\eqref{br3} are obtained easily by applying the braid-type isomorphisms \eqref{bt0}, \eqref{bt1},\eqref{bt2},\eqref{bt3} on $\cH_\bfi$. For example, take $h=(\cd,x,y,z,w,\cd)\in\cH_\bfi$ as in (2).
Then, it satisfies $-x+c_1y+z=0=-y+c_2z-w$ and then one has $X$ in \eqref{X} as
\begin{eqnarray*}
X&=& -w+(c_2x-y+w+c_2(-x+c_1y-z)_+)_+
=-w+(c_2x-y+w)_+=-w+(y-c_2z+w)_+=-w
\end{eqnarray*}
and similarly, one has $Y=-z$, $Z=-y$ and $W=-x$. 
Then \eqref{br2} has been shown.

Next, let us show that $\phi_{ij}^{(k)}(h)\in\cH_{\bfi'}$
for any $h\in\cH_\bfi$.
We shall see the case of $\phi_{ij}^{(2)}$.
The other cases will be shown similarly to the case of $\phi_{ij}^{(2)}$. 
Let $x$ be $m$-th entry in $h$, namely, 
$(i_m,i_{m+1},i_{m+2},i_{m+3})=(i,j,i,j)$.
For $\bfi=i_1\cd \stackrel{m}{i}\stackrel{\,m+1}{j}\stackrel{m+2}{i}\stackrel{\,m+3}{j}\cd i_N$, let $\bfi':=i_1\cd \stackrel{m}{j}\stackrel{\,m+1}{i}\stackrel{\,m+2}{j}\stackrel{\,m+3}{i}\cd i_N$. 
For any $h=(\cd \stackrel{m}{x},y,z,w,\cd)\in \cH_\bfi$ and
 $k\in [1,N]$ with $k^+\leq  N$ let us show that
\begin{equation}
\beta_k^{(\bfi')}(\phi_{ij}^{(2)}(h))=0.
\label{beta-h}
\end{equation}
To see all $\beta_k^{(\bfi')}(x)$'s with the variables 
$x_m,x_{m+1},x_{m+2},x_{m+3}$, it suffices to investigate the following cases:
\begin{enumerate}
\item $k=m$, that is, $\beta_m^{(\bfi')}(\phi_{ij}^{(2)}(h))=0$.
\item $k=m+1$, that is, $\beta_{m+1}^{(\bfi')}(\phi_{ij}^{(2)}(h))=0$.
\item $k=m+2$, that is, $\beta_{m+2}^{(\bfi')}(\phi_{ij}^{(2)}(h))=0$.
\item $k=m+3$, that is, $\beta_{m+3}^{(\bfi')}(\phi_{ij}^{(2)}(h))=0$.
\item $k=m^{-}$, that is, $\beta_{m^-}^{(\bfi')}(\phi_{ij}^{(2)}(h))=0$.
\item $k=(m+1)^{-}$, that is, $\beta_{(m+1)^-}^{(\bfi')}(\phi_{ij}^{(2)}(h))=0$.
\item $\beta_k^{(\bfi')}(\phi_{ij}^{(2)}(h))=0$ for $k$ such that $i_k\ne i,j$ and 
$[m,m+3]\subset [k,k^+]$.
\end{enumerate}
Since all cases are shown by a similar way, we shall only show the cases (6) and (7).\\
The case (6): We assume that there exists the index $i$ before $x=i_m=i$ in $\bfi$, 
that is, $i_{m^-}=i$.
If no such $i$ exists, there is nothing to show. Under the assumption,
note that one has $i_m=i_{m^-}=i$ for $\bfi$ and $i'_{m+1}=i'_{(m+1)^-}=i$ 
for $\bfi'=(i'_1,\cd,i'_N)$.
For $h=(h_1,\cd,h_{m^-},\cd,h_{m-1},x,y,z,w,\cd,h_N)$, set 
$h'=\phi_{ij}^{(2)}(h)=(h'_1,\cd, h'_N)$. One has 
\[
(h'_{(m+1)^-},h'_m,h'_{m+1},h'_{m+2},h'_{m+3})=(h_{m^-},w,z,y,x),\qq
h'_j=h_j \hbox{ for }j\not\in [m,m+3], 
\]
by the formula 
\eqref{br2} that is, 
$h'=\phi_{ij}^{(2)}(h)=(\cd,h_{(m+1)^-},\cd,h_{m-1},w,z,y,x,\cd)$. 
Since $\beta^{(\bfi)}_{(m+1)^-}(h)=0$, one has
\begin{eqnarray*}
&&\beta^{(\bfi)}_{m^-}(h)=h_{m^-}+\cd+a_{i\,i_{m-1}}h_{m-1}+x=0,\\
&&\beta^{(\bfi')}_{(m+1)^-}(\phi_{ij}^{(2)}(h))
=h'_{(m+1)^-}+\cd+a_{i\,i_{m-1}}h'_{m-1}+(-c_1w+z)
=h_{m^-}+\cd+a_{i\,i_{m-1}}h_{m-1}+(-c_1w+z).
\end{eqnarray*}
Thus, to show $\beta^{(\bfi')}_{(m+1)^-}(\phi_{ij}^{(2)}(h))=0$, one may get
$x=-c_1w+z$, which is immediate from the fact $x-c_1y+z=y-c_2z+w=0$ and $c_1c_2=2$.\\
\medskip
The case (7): 
%Since we consider a finite Cartan matrix 
%$A=(a_{ij})$ and then the corresponding Dynkin diagram does not include any cycle.
%Thus, for $i_k=l$ one has $a_{li}=0$ or $a_{lj}=0$. 
%So, here assume that $a_{li}\ne0$ and $a_{lj}=0$. 
Set $i_k=l$. Then one has
\begin{eqnarray*}
0&=&\beta^{(\bfi)}_k(h)=h_k+a_{li_{k+1}}h_{k+1}+\cd+a_{li}x+a_{lj}y+a_{li}z+
a_{lj}w+\cd+h_{k^+}\\
&=&h_k+a_{li_{k+1}}h_{k+1}+\cd+a_{lj}w+a_{li}z+a_{lj}y+a_{li}x+\cd+h_{k^+}=
\beta^{(\bfi')}_k(\phi_{ij}^{(2)}(h)).
\end{eqnarray*}
Therefore, we obtain that $\phi_{ij}^{(p)}(\cH_\bfi)\subset\cH_{\bfi'}$ for 
$p=0,1,2,3$. Since $\phi_{ij}^{(p)}\circ \phi_{ji}^{(p)}={\rm id}$, it yields
$\phi_{ij}^{(p)}(\cH_\bfi)=\cH_{\bfi'}$ for $p=0,1,2,3$.
\qed

%%%%%%%%%%%%%%%%%
% 
%\section{Proof: Connectedness of 
%the Cellular Crystal $\bbB_{i_1i_2\cd i_k}$}
%
%For $B(\ify)\subset \bbB_\bfi=\bbZ^N$, the condition $({\rm H}_\bfi)$ is as follows:
% {The condition (${\rm H}_\bfi$)}
%There exists the set of linear functions $\Xi_\bfi$:
%\[
%\begin{array}{l}\displaystyle
%\Xi_\bfi\subset
%\left\lan%
%\vp(x)\in (\bbQ^N)^*\,\left|\,\vp(x)=x_j-\sum_{1\leq k\leq N, \TY(k,m,)
%\leq N}c_k\beta_k(x)\,\,(c_k\in\bbZ,\,\,1\leq j\leq N)\right.
%% 
%\right\ran_{+}\\
%
%s.t. \,\, B(\ify)=\{x\in \bbZ^N=\bbB_\bfi\,|\,
%\vp(x)\geq0\,\,(\vp\in\Xi_\bfi)\}
%\end{array},
%\]
%where $\lan S\ran_+$ is the set of all positive linear combinations
%of the set $S$.
% 
%\begin{thm}
%If the condition $({\rm H}_\bfi)$ holds for $\bfi$, then 
%the cellular crystal $\bbB_\bfi=B_{i_1}\ot B_{i_2}\ot\cd\ot B_{i_N}$ is connected.
%\end{thm}
%This theorem is shown by the following facts:
% 
%%%%%%%%%%%%%%%%%%
 
%By the fact that {$B(\ify)$ is connected} 
In \cite[Sect.8]{Kana-N}, we have shown the following statements 
under the condition "${\bf H}_\bfi$", where we omit 
the explicit form of ${\bf  H}_\bfi$ since we do not need it here.
But, we succeed in showing the following proposition 
without the condition ${\bf H}_\bfi$ 
since in \cite{Kana-N} we have shown that 
there exists a specific reduced longest word 
$\bfi_0$ satisfying the condition ${\bf H}_{\bfi_0}$ for each simple Lie algebra $\ge$ and we 
got Lemma \ref{310}.
\begin{pro}\label{pro-bi}
Let $\bfi=i_1i_2\cd i_N$ be an arbitrary reduced longest word.
Here if the crystal $B(\ify)$ is realized in $\bbB_\bfi$ as in \ref{subsec-half}, 
we shall denote 
it by $B(\ify)_\bfi$ to emphasize the word $\bfi$. 
For $h\in\cH_\bfi$, define 
\[
B^h(\ify)_\bfi:=\{x+h\in \bbZ^N(=\bbB_\bfi)\,|\,x\in B(\ify)_\bfi\}\subset \bbB_\bfi.
\]
\begin{enumerate}
\item
For any $x+h\in B^h(\ify)_\bfi$ and $i\in I$, we obtain
\begin{equation}
\eit(x+h)=\eit(x)+h,\qq
\fit(x+h)=\fit(x)+h.
\label{eifih}
\end{equation}
%and then  $B^h(\ify)$ is connected.
\item
For any $h\in \cH_\bfi$, we have\,\,
$B(\ify)_\bfi\cap B^h(\ify)_\bfi\ne\emptyset$.
\item
\[
\bbB_\bfi=\bigcup_{h\in\cH_\bfi}B^h(\ify)_\bfi
\]
\end{enumerate}
\end{pro}
\begin{rem}
In the setting of the half-potential method in \cite{Kana-N}, as mentioned in
Remark \ref{trop-real}, the crystal 
$B(\ify)$ is realized as a subset of $\bbB_\bfi$ and it is supposed that 
$\eit x=0$ if $\eit x\not\in (\wtil \bbB^-_{w_0})_{\Phi^{(+)},\Theta_\bfi}\cong B(\ify)$. 
At the statement (2), since 
$x\in B(\ify)_\bfi$ is considered as an element of $\bbB_\bfi$, 
$\eit x$ is also considered as an element in $\bbB_\bfi$. That is, 
even if $\eit x\not\in B(\ify)$, we consider that $\eit x\in\bbB_\bfi$ and then
it never vanishes. 
\end{rem}
{\sl Proof.} 
First, note that in \cite{Kana-N}, it has been shown the statement (1)-(3) 
for some specific reduced longest word $\bfi_0$ satisfying 
the condition ${\bf H}_{\bfi_0}$. 

\medskip
\nd
(1) One can show \eqref{eifih} easily since $\beta^{(\bfi)}_k(h)=0$ for any 
$k\in[1,N]$ and $h\in \cH_\bfi$ and then one has 
$\TY(\wtil m,i,e)(x+h)=\TY(\wtil m,i,e)(x)$ and 
$\TY(\wtil m,i,f)(x+h)=\TY(\wtil m,i,f)(x)$. 

\medskip
\nd
(2) 
For any reduced longest word $\bfi$, 
there exists a composition of braid-moves $\xi$ on $\bfi_0$ such that 
$\xi(\bfi_0)=\bfi$. Let $\Xi$ be the composition of braid-type isomorphisms 
associated with $\xi$ such that $\Xi(B_{\bfi_0})=B_\bfi$.
Since $\Xi$ is an isomorphism of crystals and then commutes with the actions of 
$\eit$ and $\fit$ and 
$B(\ify)_{\bfi_0}\cap B^h(\ify)_{\bfi_0}\ne\emptyset$ for any $h\in \cH_{\bfi_0}$,
one gets that 
$\Xi(B^h(\ify)_{\bfi_0})=B^{\Xi(h)}(\ify)_\bfi$ for any $h\in \cH_{\bfi_0}$, 
which implies that  
$B(\ify)_\bfi\cap B^H(\ify)_\bfi\ne\emptyset$
for any $H\in \cH_\bfi$ since $\Xi(\cH_{\bfi_0})=\cH_\bfi$ by Lemma \ref{310}.

\medskip
\nd
(3) Let $\xi$ and $\Xi$ be as in (2). 
Then, since $\Xi(\cH_{\bfi_0})=\cH_\bfi$ and 
$\Xi(B^h(\ify)_{\bfi_0})=B^{\Xi(h)}(\ify)_\bfi$ for any $h\in \cH_{\bfi_0}$,
we obtain 
\[
\qq\qq\qq\bbB_\bfi=\Xi(\bbB_{\bfi_0})
=\bigcup_{h\in\cH_{\bfi_0}}\Xi(B^h(\ify)_{\bfi_0}) 
=\bigcup_{h\in\cH_{\bfi_0}}B^{\Xi(h)}(\ify)_{\bfi}
=\bigcup_{H\in\cH_{\bfi}}B^H(\ify)_{\bfi}.\qq\qq\qq\qq\qed
\]
 It is immediate from  this proposition that one has the following theorem:
\begin{thm}[\cite{Kana-N}]
For any simple Lie algebra $\ge$ and any reduced word $i_1i_2\cd i_k$, 
the 
cellular crystal $\bbB_{i_1i_2\cd i_k}=B_{i_1}\ot B_{i_2}\ot\cd\ot B_{i_k}$
is connected as a crystal graph.
\end{thm}

\section{Quiver Hecke Algebra and its modules}\label{QHA}
\subsection{Definition of Quiver Hecke Algebra}\label{defQHA}
We review basics of the quiver Hecke algebras (see e.g.,\cite{K-L,Loc,jams,Rou}).
For a finite index set $I$ and a field $\bf k$, 
let $({\mathscr Q}_{i,j}(u,v))_{i,j\in I}\in{\bf k}[u,v]$ be polynomials satisfying:
\begin{enumerate}
\item ${\mathscr Q}_{i,j}(u,v)={\mathscr Q}_{j,i}(v,u)$ for any $i,j\in I$.
\item ${\mathscr Q}_{i,j}(u,v)$ is in the form:
\[
{\mathscr Q}_{i,j}(u,v)=\begin{cases}\displaystyle
\sum_{a(\al_i,\al_i)+b(\al_j,\al_j)=-2(\al_i,\al_j)}
t_{i,j;a,b}u^av^b&\hbox{ if }i\ne j,\\
0&\hbox{ if }i=j,
\end{cases}
\]
where $t_{i,j;-a_{ij},0}\in{\bf k}^\times$.
\end{enumerate}
For $\beta=\sum_im_i\al_i\in Q_+$ with $|\beta|:=\sum_im_i=m$, set
$I^\beta:=\{\nu=(\nu_1,\cd,\nu_m)\in I^m\mid \sum_{k=1}^m\al_{\nu_k}=\beta\}$.
\begin{df}
For $\beta\in Q_+$, the {\it quiver Hecke algebra $R(\beta)$} associated with a Cartan matrix  $A$ and 
polynomials ${\mathscr Q}_{i,j}(u,v)$ is the $\bfk$-algebra generated by
\[
\{e(\nu)|\nu\in I^\beta\},\q
\{x_k|1\leq k\leq n\},\q
\{\tau_i|1\leq i\leq n-1\}
\]
with the following relations:
\begin{eqnarray*}
&&e(\nu)e(\nu')=
\del_{\nu,\nu'}e(\nu),\q
\sum_{\nu\in I^\beta}e(\nu)=1, \q
e(\nu)x_k=x_ke(\nu), \q
x_kx_l=x_lx_k,\\
&&\tau_le(\nu)=e(s_l(\nu))\tau_l,\q
\tau_k\tau_l=\tau_l\tau_k\,\,\hbox{ if }|k-l|>1,\\
&& \tau_k^2e(\nu)={\mathscr Q}_{\nu_k,\nu_{k+1}}(x_k,x_{k+1})e(\nu),\\
&&(\tau_kx_l-x_{s_k(l)}\tau_k)e(\nu)=\begin{cases}
-e(\nu)&\hbox{ if }l=k, \,\,\nu_k=\nu_{k+1},\\
e(\nu)&\hbox{ if }l=k+1, \,\,\nu_k=\nu_{k+1},\\
0&\hbox{otherwise},
\end{cases}\\
&&
(\tau_{k+1}\tau_k\tau_{k+1}-\tau_k\tau_{k+1}\tau_k)e(\nu)=
\begin{cases}
\ovl{\mathscr Q}_{\nu_k,\nu_{k+1}}(x_k,x_{k+1},x_{k+2})e(\nu)
&\hbox{ if } \nu_k=\nu_{k+2},\\
0&\hbox{ otherwise},
\end{cases}
\end{eqnarray*}
where $\ovl{\mathscr Q}_{i,j}(u,v,w)=\frac{{\mathscr Q}_{i,j}(u,v)-{\mathscr Q}_{i,j}(w,v)}{u-w}
\in\bfk[u,v,w]$.
\end{df}
%and set $ R:=\bigoplus_{\beta\in Q_+}R(\beta)$.
\begin{enumerate}
\item
The relations above are homogeneous if we define 
\[
{\rm deg}(e(\nu))=0,\q
{\rm deg}(x_ke(\nu))=(\al_{\nu_k},\al_{\nu_k}),\q
{\rm deg}(\tau_l e(\nu))=-(\al_{\nu_l},\al_{\nu_{l+1}}).
\]
Thus, $R(\beta)$ becomes a $\bbZ$-graded algebra. Here we define the {\it weight} of 
$R(\beta)$-module $M$ as ${\rm wt}(M)=-\beta$.
\item Let $M=\bigoplus_{k\in\bbZ}M_k$ be a $\bbZ$-graded $R(\beta)$-module. Define
a {\it grading shift functor} $q$ on  
the category of graded $R(\beta)$-modules $R(\beta)$-Mod by 
\[
qM:=\bigoplus_{k\in\bbZ}(qM)_k,\q\hbox{where }
(qM)_k=M_{k-1}.
\]
\item
For $M,\,N\in R(\beta)$-Mod, let
${\rm Hom}_{R(\beta)}(M,N)$ be the space of degree preserving morphisms and define 
$\textsc{Hom}_{R(\beta)}(M,N):=\bigoplus_{k\in \bbZ}{\rm Hom}_{R(\beta)}(q^kM,N)$, 
which is a space of morphisms up to grading shift. We define
${\rm deg}(f)=k$ for $f\in {\rm Hom}_{R(\beta)}(q^kM,N)$. 
\item
Let $\psi$ be the anti-automorphism of $R(\beta)$ preserving all generators. 
For $M\in R(\beta)$-Mod, define $M^*:=\textsc{Hom}_\bfk(M,\bfk)$ with the $R(\beta)$-
module structure by $(r\cdot f)(u):=f(\psi(r)u)$ for 
$r\in R(\beta)$, $u\in M$ and $f\in M^*$, which is called a {\it dual module} of $M$.
In particular, if $M\cong M^*$ we call $M$ is {\it self-dual}.
\item
For $\beta,\gamma\in Q_+$, set $e(\beta,\gamma)
=\sum_{\nu\in I^\beta,\nu'\in I^\gamma}e(\nu,\nu')$. 
We define an injective homomorphism
$\xi_{\beta,\gamma}:R(\beta)\ot R(\gamma)\to 
e(\beta,\gamma)R(\beta+\gamma)e(\beta,\gamma)$ by 
$\xi(\beta,\gamma)(e(\nu)\ot e(\nu'))=e(\nu,\nu')$, 
$\xi(\beta,\gamma)(x_k e(\beta)\ot1)=x_ke(\beta,\gamma)$,
$\xi(\beta,\gamma)(1\ot x_ke(\gamma))=x_{k+|\beta|}e(\beta,\gamma)$,
$\xi(\beta,\gamma)(\tau_ke(\beta)\ot1)=\tau_ke(\beta,\gamma)$,
$\xi(\beta,\gamma)(1\ot\tau_ke(\gamma))=\tau_{k+|\beta|}e(\beta,\gamma)$.
\item
For $M\in R(\beta)$-Mod and $N\in R(\gamma)$-Mod, define 
the {\it convolution product} $\circ$ by 
\[
M\circ N:=R(\beta+\gamma)e(\beta,\gamma)\ot_{R(\beta)\ot R(\gamma)}
(M\ot N)
\]
For simple $M\in R(\beta)$-Mod and simple $N\in R(\gamma)$-Mod, we say $M$ and $N$ 
{\it strongly commutes} if $M\circ N$ is simple and $M$ is {\it real} if 
$M\circ M$ is simple.
\item
For $M\in R(\beta)$-Mod and $N\in R(\gamma)$-Mod, denote by 
$ M\nabla N:={\rm hd}(M\circ N)$ the head of $M\circ N$ and 
$ M\Del N:={\rm soc}(M\circ N)$ the socle of $M\circ N$, 
where the head of module $M$ is the quotient by its radical and the socle of module $M$ is the summation of all simple submodules.
\end{enumerate}

%%%%%%%%%%%%%%%%%%%%%%%%%%%%%%%%%%%%%%%%%%%
 %[allowframebreaks]
\subsection{Categorification of quantum coordinate ring $\aq$}

Let $R(\beta)$-{\rm gmod} be the full subcategory of $R(\beta)$-Mod whose 
objects are  
finite-dimensional graded $R(\beta)$-modules and set
$R$-gmod$=\bigoplus_{\beta\in Q_+}R(\beta)$-gmod.
% and $R(\beta)$-{\rm gproj} the full-subcategory of $R(\beta)$-Mod whose objects are 
%finitely generated projective $R(\beta)$-modules.
Define the functors 
\[
%{ E_i}&:&R(\beta)\hbox{-gproj}\to R(\beta-\al_i)\hbox{-gproj}, \q
E_i:R(\beta)\hbox{-gmod}\to R(\beta-\al_i)\hbox{-gmod}, \qq
%{ F_i}&:&R(\beta)\hbox{-gproj}\to R(\beta+\al_i)\hbox{-gproj}, \q
F_i:R(\beta)\hbox{-gmod}\to R(\beta+\al_i)\hbox{-gmod }, 
\]
by $E_i(M):=e(\al_i,\beta-\al_i)M,\q
F_i(M)=L(i)\circ M$,
where $ e(\al_i,\beta-\al_i):=
\sum_{\nu\in I^\beta, \nu_1=i}e(\nu)$ and 
$L(i):=R(\al_i)/R(\al_i)x_1$ is a 1-dimensional simple $R(\al_i)$-module. 
Let $\cK(R\hbox{\rm -gmod})$ be the Grothendieck ring of $R$-gmod and 
then $\cK(R\hbox{\rm -gmod})$ becomes a $\bbZ[q,q^{-1}]$-algebra
with the multiplication induced by
the convolution product and $\bbZ[q,q^{-1}]$-action induced by the grading shift functor $q$. Here, one obtain the following:
\begin{thm}[\cite{K-L,Rou}]
As a $\bbZ[q,q^{-1}]$-algebra there exists an isomorphism
\[
\cK(R\hbox{\rm -gmod})\cong \aq_{\bbZ[q,q^{-1}]}.
\]
\end{thm}

%%%%%%%%%%%%%%%%%%%
%%%%%%%%%%%%%%%%%%%%%%%%%%%%%%%%%%%%%%%%%%%
\subsection{Categorification of the crystal $B(\ify)$ by Lauda and Vazirani \cite{L-V}}
The following lemma is given in \cite{K-L}:
\begin{lem}[\cite{K-L}]\label{simple}
For any simple $R(\beta)$-module $M$, ${\rm soc}(E_i M)$, ${\rm hd}(E_i M)$
 and ${\rm hd}(F_i M)$ are all simple modules. Here we also have that
${\rm soc}(E_i M)\cong{\rm hd}(E_i M)$ up to grading shift.
\end{lem}
For $M\in R(\beta)$-gmod, define
\begin{eqnarray}
&&\wt(M)=-\beta,\q
\vep_i(M)=\max\{n\in\bbZ\,|\,E_i^nM\ne0\},\q 
\vp_i(M)=\vep_i(M)+\lan h_i,\wt(M)\ran,\\
&&\Eit M:= q_i^{1-\vep_i(M)}{\rm soc}(E_i M)\cong q_i^{\vep_i(M)-1}{\rm hd}(E_iM),\qq
\Fit M:=q_i^{\vep_i(M)}{\rm hd}(F_iM).\label{eitfit}
\end{eqnarray}
Set 
$\bbB(\rgmod):=\{S\,|\,S\hbox{ is a self-dual simple module in }R\hbox{-gmod}\}$.
Then, it follows from Lemma \ref{simple}
that $\wtil E_i$ and $\wtil F_i$ are well-defined on $\bbB(\rgmod)$.
\begin{thm}[\cite{L-V}]\label{LV}
The 6-tuple $(\bbB(\rgmod),\{\Eit\},\{\Fit\},\wt,\{\vep_i\},\{\vp_i\})_{i\in I}$ holds 
a crystal structure and there exists  the following 
isomorphism of crystals:
\[
\Psi:\bbB(\rgmod)\q\mapright{\sim}\q B(\ify).
\]
\end{thm}
\begin{rem}
Note that Lauda and Vasirani showed this theorem 
under more general setting that $\ge$ is arbitrary symmetrizable Kac-Moody Lie algebra.
Here we assume that $\ge$ is a simple Lie algebra.
The definition of $\Eit$ and $\Fit$ in \eqref{eitfit} differs from the one in \cite{L-V}, which follows the one in \cite{jams}.
\end{rem}

%%%%%%%%%%%%%%%%%%%%%%%%%%%%%%%%%%%%%%%%%%%%%%%%
\section{Localization of monoidal category}
%%%%%%%%%%%%%%%%%%%%%%%%%%%%%%%%%%%%%%%%%%
In this section, we shall review briefly the basics on localization of monoidal category following \cite{Loc}. 
\subsection{Braiders and Real Commuting Family}
Let $\Lm$ be $\bbZ$-lattice and $\cT=\oplus_{\lm\in \Lm}
\cT_\lm$ be a $\bfk$-linear $\Lm$-graded monoidal category with a data consisting
of a bifunctor $\ot :\cT_\lm\times \cT_\mu\to \cT_{\lm+\mu}$, an isomorphism 
$a(X,Y,Z):(X\ot Y)\ot Z\mapright{\sim}X\ot(Y\ot Z)$ satisfying  
$a(X,Y,Z\ot W)\circ a(X\ot Y,Z,W)={\rm id}_X\ot a(Y,Z,W)\circ a(X,Y\ot Z,W)\circ a(X,Y,Z)\ot {\rm id}_W$ and an object ${\bf 1}\in\cT_0$
endowed with an isomorphism $\epsilon:{\bf 1}\ot{\bf 1}\mapright{\sim}{\bf 1}$ 
such that the functor $X\mapsto X\ot {\bf 1}$ 
and $X\mapsto {\bf1}\ot X$ are fully-faithful. 
%(Later $\Lm$ will be the root lattice $Q$)
\begin{df}[\cite{Loc}]
Let $q$ be the grading shift functor on $\cT$. A {\it graded braider} is a triple
$(C,R_C,\phi)$, where $C\in \cT$, $\bbZ$-linear map $\phi:\Lm\to\bbZ$ and a 
morphism:
\[
R_C:C\ot X\to q^{\phi(\lm)}X\ot C
\q(X\in\cT_\lm),
\]
satisfying the following commutative diagram:
\[
\xymatrix{
\linethickness{30pt}
C\ot X\ot Y\ar@{->}^{R_C(X)\ot Y}[r]
\ar@{->}_{R_C(X\ot Y)}[dr]& q^{\phi(\lm)}\ot X\ot C\ot Y
\ar@{->}^{X\ot R_C(Y)}[d]\\
&q^{\phi(\lm+\mu)}(X\ot Y)\ot C
}\qq
(X\in \cT_\lm,\,\,Y\in\cT_\mu)
\]
and being functorial, that is, for any $X,Y\in \cT$ and $f\in{\rm Hom}_\cT(X,Y)$
it  satisfies the following  commutative diagram:
\[
\xymatrix{
\linethickness{30pt}
C\ot X\ar@{->}^{\id\ot f}[r]
\ar@{->}_{R_C(X)}[d]
&C\ot Y\ar@{->}^{R_C(Y)}[d] \\
X\ot C \ar@{->}^{f\ot\id}[r] & Y\ot C
}
\]
\end{df}

\begin{df}[\cite{Loc}]\label{real-com}
Let $I$ be an  index set and $(C_i,R_{C_i},\phi_i)_{i\in I}$ a 
family of graded braiders in $\cT$. We say that 
$(C_i,R_{C_i},\phi_i)_{i\in I}$ is a
{\it real commuting family of graded braiders} in $\cT$ if
\begin{enumerate}
\item
$C_i\in \cT_{\lm_i}$ for some $ \lm_i\in\Lm$, 
and $ \phi_i(\lm_i)=0$, 
$ \phi_i(\lm_j)+\phi_j(\lm_i)=0$ for any $i,j\in I$.
\item
$R_{C_i}(C_i)\in\bfk^\times {\rm id}_{C_i\ot C_i}$  for any $i\in I$.
\item
$R_{C_i}(C_j)\ot R_{C_j}(C_i)
\in\bfk^\times {\rm id}_{C_i\ot C_j}$ 
for any $i,j\in I$.
\end{enumerate}
\end{df}
Note that $R_{C_i}$'s satisfy so-called "Yang-Baxter equation", such as, 
\[
R_{C_i}(C_j)\circ R_{C_i}(C_k)\circ R_{C_j}(C_k)
=R_{C_j}(C_k)\circ R_{C_i}(C_k)\circ R_{C_i}(C_j) \,\,\hbox{ on }\,\,
C_i\circ C_j\circ C_k.
\]
For a finite index set $I$,  set 
$\Gamma:=\oplus_{i\in I}\bbZ e_i$ and 
$\Gamma_+:=\oplus_{i\in I}\bbZ_{\geq0} e_i$.
%(Later $\Gamma$ will be the weight lattice and $\Gamma_+$ be the set of 
%dominant weights.)
\begin{lem}[\cite{Loc}]
Suppose that we have a real commuting family of graded braiders 
$(C_i,R_{C_i},\phi_i)_{i\in I}$.
We can choose a bilinear map $H:\Gamma\times \Gamma\to \bbZ$ such that 
$\phi_i(\lm_j)=H(e_i,e_j)-H(e_j,e_i)$ and there exist
\begin{enumerate}
\item an object $C^\al$ for any $\al\in \Gamma_+$.
\item
an isomorphism $\xi_{\al,\beta}:C^\al\ot  C^\beta \mapright{\sim}
q^{H(\al,\beta)}C^{\al+\beta}$ for any $\al,\beta\in\Gamma_+$
\end{enumerate}
such that $C^0=1$ and $C^{e_i}=C_i$.
% and satisfying some commutative diagrams.
\end{lem}

%%%%%%%%%%%%%%%%%%%%%%%%%%%%%%%%%%%%%%%%%%%
\subsection{Localization}

Let $\cT$ and $(C_i,R_{C_i},\phi_i)_{i\in I}$
be  as above and 
$\{C^\al\}_{\al\in\Gamma_+}$ objects as in the previous lemma.
We define a partial order $\preceq$ on $\Gamma$ by 
\[
\al\preceq \beta\Longleftrightarrow \beta-\al\in\Gamma_+
\]
For $\al_1,\al_2,\cd\in\Gamma$, define
\[
\cD_{\al_1,\al_2,\cd}:=\{\del\in\Gamma\,|\,
\al_j+\del\in \Gamma_+\,\,\hbox{ for any } j=1,2,\cd\}.
\]
For $X\in \cT_\lm$, $Y\in \cT_\mu$ and $\del\in \cD_{\al,\beta}$,
set 
\[
H_\del((X,\al),(Y,\beta))
:={\rm Hom}_{\cT}(C^{\del+\al}\ot X,
q^{P(\al,\beta,\del,\mu)}Y\ot C^{\beta+\del}),
\]
where a $\bbZ$-valued function $P(\al,\beta,\del,\mu):=
H(\del,\beta-\al)+\phi(\del+\beta,\mu)$ and the map $\phi:\Gamma\times \Lm\to\bbZ$
is defined by $\phi(\al,L(\beta))=H(\al,\beta)-H(\beta,\al)$ and 
$L:\Gamma\to\Lm$ is defined by $L(e_i)=\lm_i$ (\cite{Loc}).

\begin{lem}[\cite{Loc}]
For $\del\preceq \del'$ there exists the map 
\[
\zeta_{\del,\del'}:H_\del((X,\al),(Y,\beta))
\to
H_{\del'}((X,\al),(Y,\beta))
\]
satisfying
\[
\zeta_{\del,\del'}\circ \zeta_{\del',\del''}=\zeta_{\del,\del''}
\hbox{for }\del\preceq\del'\preceq \del''.
\]
Therefore, we find that 
$\{H_\del((X,\al),(Y,\beta))\}_{\del\in \cD_{\al,\beta}}$ 
becomes an inductive system.
\end{lem}
\begin{df}[Localization \cite{Loc}]
We define the category $\wtil\cT$ by 
\begin{eqnarray*}
&&{\rm Ob}(\wtil\cT):={\rm Ob}(\cT)\times\Gamma,\\
&&\Hom_{\wtil\cT}((X,\al),(Y,\beta)):=\varinjlim_{
\begin{array}{c}\scriptstyle\tiny \del\in\cD(\al,\beta),\\ 
\scriptstyle\lm+L(\al)=\mu+L(\beta)
\end{array}}
H_\del((X,\al),(Y,\beta)),
\end{eqnarray*}
where $X\in\cT_\lm$, $Y\in \cT_\mu$ and the function 
$L:\Gamma\to\Lm$ ($e_i\mapsto \lm_i$) is as above.
We call this $\wtil\cT$ a {\it localization} of $\cT$ by 
$(C_i,R_{C_i},\phi_i)_{i\in I}$ and 
denote it by $\cT[C_i^{\ot -1}\,|\,i\in I]$ when we emphasize 
$\{C_i\mid i\in I\}$.
\end{df}
%%%%%%%%%%%%%%%%%%%%%%%%%%%%%%%%%%%%%%%%%%
\begin{thm}[\cite{Loc}]\label{til}
$\wtil\cT$ becomes a monoidal category. Moreover, there exists a monoidal 
functor $\Upsilon:\cT\to \wtil\cT$ such that 
\begin{enumerate}
\item 
$\Upsilon(C_i)$ is {\it invertible} in $\wtil\cT$ for any $i\in I$, 
namely, the functors
 $X\mapsto X\ot\Upsilon(C_i)$ and $X\mapsto \Upsilon(C_i)\ot X$ 
are equivalence of categories.
\item
For any $i\in I$ and $X\in \cT$, $\Upsilon(R_{C_i}(X)):
\Upsilon(C_i\ot X)\to\Upsilon(X\ot C_i)$ is an isomorphism.
\item
The functor $\Upsilon$ holds the following universality:
If there exists another monoidal category $\cT'$ and a monoidal fucntor 
$\Upsilon':\cT\to\cT'$ satisfying the above statements (1) and (2), then
there exists a monoidal functor $F:\wtil\cT\to\cT'$ (unique up to iso.)
such that $\Upsilon'=F\circ \Upsilon$. 
\end{enumerate}
\end{thm}

\begin{pro}[\cite{Loc}]\label{til-pro}
Under the setting above, we obtain
\begin{enumerate}
\item
$(X,\al+\beta)\cong q^{-H(\beta,\al)}(C^\al\ot X,\beta)$, 
$(1,\beta)\ot (1,-\beta)\cong q^{-H(\beta,\beta)}(1,0)$ for
$\al\in\Gamma_+$, $\beta\in\Gamma$ and $X\in \wtil\cT$.
\item
If $\cT$ is an abelian category, then so is $\wtil\cT$.
\item
The functors $\Upsilon:\cT\to\wtil\cT$ is exact.
\item
If the functor $-\ot Y$ and $Y\ot -$ are exact for any $Y$ in $\cT$, 
then the functors 
$\wtil\cT\to\wtil\cT$ ($X\mapsto X\ot Y$ (resp. $X\to Y\ot X$)) are 
exact for any $Y$ in $\wtil\cT$.

\end{enumerate}

\end{pro}
 
%%%%%%%%%%%%%%%%%%%%%%%%%%%%%%%%%%%%%%%%%%%
% %[allowframebreaks]
%\section{Cyclotomic Quiver Hecke Algebra s}
%For a dominat weight $\lm\in P_+$, set 
%$a_\Lm:=\{a_{\Lm,i}(t_i)=t_i^{\lan h_i,\Lm\ran}\}$ .
%Define the {\bf cyclotomic quiver Hecke algebra} $R^{a_\Lm}(\lm)$
%($\lm\in P$) by 
%\[
%R^{a_\Lm}(\lm):=\frac{R(\beta)}{\sum_iR(\beta)
%a_{\Lm,i}(x_n e(\beta-\al_i,\al_i)R(\beta)}
%\q(\beta=\Lm-\lm)
%\]
%We define the functors $F_i^{a_\Lm}$, $E_i^{a_\Lm}$, 
%$F_i^{a_\Lm\,(n)}$, $E_i^{a_\Lm\,(n)}$ by 
%\begin{eqnarray*}
%&&F_i^{a_\Lm}M:=R^{a_\Lm}(\lm-\al_i)e(\al_i,\beta)\ot_{R^{a_\Lm}(\lm)}M,\\
%&&E_i^{a_\Lm}M:=e(\al_i,\beta-\al_i)\ot_{R^{a_\Lm}(\lm)}M,\\
%&&
%F_i^{a_\Lm(n)}M:=\Hom_{R(n\al_i)}(P(i^n),(F_i^{a_\Lm})^nM),\\
%&&E_i^{a_\Lm(n)}M:=\Hom_{R(n\al_i)}(P(i^n),(E_i^{a_\Lm})^nM)
%\end{eqnarray*}
%
%
%
% 
%%%%%%%%%%%%%  Sect.6 %%%%%%%%%%%%%%%%%%%%%%%%%%%%%
\section{Localization of the category $R$-{\rm gmod}}
\subsection{Determinantial Modules}
Here we just go back to the setting as in Sect.\ref{QHA}.
Let  $L(i^n):=q_i^{\frac{n(n-1)}{2}}L(i)^{\circ n}$ 
be a simple $R(n\al_i)$-module satisfying  
${\rm qdim}(L(i^n))=[n]_i!:=\prod_{k=1}^n\frac{q_i^k-q_i^{-k}}{q_i-q_i^{-1}}$
\,\,($q_i:=q^{\frac{(\al_i,\al_i)}{2}}$).
\begin{df}[\cite{Loc,jams}]
For $M\in R$-gmod, define 
\[
\wtil F_i^n(M):=L(i^n)\nabla M.
\]
For a Weyl group element $w$, let $s_{i_1}\cd s_{i_l}$ 
be its reduced expression. For a dominant weight $\Lm\in P_+$, 
set 
\[
m_k:=\lan h_{i_k},s_{i_{k+1}}\cd s_{i_l}\Lm\ran\qq (k=1,\cd,l).
\] 
We define the
{\it determinantial module} associated with $w$ and $\Lm$ by
\[
{\bf M}(w\Lm,\Lm):=\wtil F_{i_1}^{m_1}\cd
\wtil F_{i_l}^{m_l}{\bf 1},
\]
where ${\bf 1}$ is a trivial $R(0)$-module.
\end{df}
Note that in general, one can define determinantial modules ${\bf M}(w\Lm,u\Lm)$
($w,u\in W$) which corresponds to the generalized minor $\Del_{w\Lm,u\Lm}$.

Now, let us see some similarity between the family of 
determinantial modules $\{{\bf M}(w_0\Lm,\Lm)\}_{\Lm\in P_+}$ and the subspace $\cH_\bfi$.
As has seen above that for a reduced longest word $\bfi=i_1\cd i_N$,
the subspace 
$\cH_\bfi\subset \bbB_\bfi$ is presented by
\[
\cH_\bfi=\bigoplus_{i\in I}\bbZ\bfh_i,\q
\bfh_i=((\TY(h,k,i):=\lan h_{i_k},s_{i_{k+1}}\cd s_{i_N}\Lm_i\ran)_{k=1,\cd,N}.
\]
Furthermore, we also get 
\begin{pro}\label{pro-f0}
For any reduced longest word $\bfi=i_1i_2\cd i_N$ and $\Lm\in P_+$, set 
\[
m_k:=\lan h_{i_k},s_{i_{k+1}}s_{i_{k+2}}\cd s_{i_l}\Lm\ran\q
(k=1,2,\cd, N)\q\hbox{ and }\q
{\bf h}_\Lm:=(m_1,\cd,m_N).
\]
Then we obtain
\begin{eqnarray*}
%\label{f-bfi}
{\bf h}_\Lm=\til f_{i_1}^{m_1}\til f_{i_2}^{m_2}\cd \til f_{i_N}^{m_N}
((0)_{i_1}\ot (0)_{i_2}\ot\cd\ot (0)_{i_N})
&=&\til f_{i_1}^{m_1}(0)_{i_1}\ot \til f_{i_2}^{m_2}(0)_{i_2}\ot\cd\ot 
\til f_{i_N}^{m_N}(0)_{i_N}\in\cH_\bfi,
\end{eqnarray*}
where note that for $\Lm=\sum_ia_i\Lm_i$, one has ${\bf h}_\Lm=
\sum_ia_i{\bf h}_{\Lm_i}$.
%Then in this case we obtain 
%\[
%\til f_{i_1}^{m_1}\til f_{i_2}^{m_2}\cd \til f_{i_N}^{m_N}
%((0)_{i_1}\ot (0)_{i_2}\ot\cd\ot (0)_{i_N})
%=\bfh_i.
%\]
\end{pro}
{\sl Proof.}
We shall show the following by descending induction on $k\in[1,N]$:
\begin{equation}
\til f_{i_k}^{m_k}\cd \til f_{i_N}^{m_N}
((0)_{i_1}\ot (0)_{i_2}\ot\cd\ot (0)_{i_N})
=(0)_{i_1}\ot\cd\ot (0)_{i_{k-1}}\ot \til f_{i_k}^{m_k}(0)_{i_k}\ot\cd\ot 
\til f_{i_N}^{m_N}(0)_{i_N}=:M_k.
\label{f-m}
\end{equation}
By Proposition \ref{basis-h}, we find that $(m_1,\cd,m_N)\in\cH_\bfi$
and for $\{k=k_0<k_1<\cd<k_L\}:=\{l\mid i_l=i_k,l\geq k\}$ 
\[
\sigma_{k}(M_{k+1})=\sigma_{k_1}(M_{k+1})=\cd =\sigma_{k_L}(M_{k+1}).
\]
Then, we get 
\[
0=\sigma_k(M_{k+1})>-m_k=\sigma_{k_1}(M_{k+1})=\cd
=\sigma_{k_L}(M_{k+1}),
\]
which implies
\[
\til f_{i_k}^{m_k}M_{k+1}
=\til f_{i_k}^{m_k}\cd \til f_{i_N}^{m_N}
((0)_{i_1}\ot (0)_{i_2}\ot\cd\ot (0)_{i_N})
=(0)_{i_1}\ot\cd\ot (0)_{i_{k-1}}\ot \til f_{i_k}^{m_k}(0)_{i_k}\ot\cd\ot 
\til f_{i_N}^{m_N}(0)_{i_N}=M_k.
\]
Thus, we finish the proof.\qed

By this proposition, one observes that there would exist a certain 
correspondence 
\begin{equation}
{\bf M}(w_0\Lm,\Lm)=\wtil F_{i_1}^{m_1}\cd
\wtil F_{i_l}^{m_l}{\bf 1}\q \longleftrightarrow \q{\bf h}_\Lm=\til f_{i_1}^{m_1}\cd \til f_{i_N}^{m_N}((0)_{i_1}\ot (0)_{i_2}\ot\cd\ot (0)_{i_N}).
\label{corr}
\end{equation}

%%%%%%%%%%%%%%%%%%%%%%%%%%%%%%%%%%%%%%%%%%
 %[allowframebreaks]
%\subsection{Affinizations}
\begin{df}[\cite{Loc}]
For $\beta\in Q_+$, define a central element in $R(\beta)$ by\\ 
${\mathfrak p}_i:=\sum_{\nu\in I^\beta}
\left(\prod_{a\in\{1,2,\cd,{\rm ht}(\beta)\},\nu_a=i}
x_a\right)e(\nu)\in R(\beta)$.
For a simple $M\in R(\beta)$-gmod, define an {\it affinization} 
$\what M$ of $M$ with degree $d$:
\begin{enumerate}
\item
There is an endomorphism $z:\what M\to\what M$ of degree $d>0$ such that 
$\what M$ is finitely generated free module of ${\bf k}[z]$ 
and $\what M/z\what M\cong M$.
\item
${\mathfrak p}_i\what M\ne0$ for any $i\in I$.
\end{enumerate}
\end{df}
\begin{thm}[{\cite[Theorem 3.26]{Loc}}]\label{aff}
For any $\Lm\in P_+$ and $w\in W$, the determinantial module ${\bf M}(w\Lm,\Lm)$
is a real simple module and admits an affinization $\what{\bf M}(w\Lm,\Lm)$.
\end{thm}

Note that indeed, if $\ge$ is simply-laced, then the affinization $\what M$ 
always exists for any simple $M\in R(\beta)$-gmod as (\cite{KP}),
\[
\what M={\bf k}[z]\ot_ {\bf k} M.
%\wtil F_{i_1}^{m_1}\cd
%\wtil F_{i_l}^{m_l}
\]
%\begin{pro}[\cite{Loc}???]
%For any $w\in W$ and $\Lm\in P_+$, 
%${\bf M}(w\Lm,\Lm)$ is a real simple module and admits an affinization, indeed, 
%\[
%{\bf M}(w\Lm,\Lm):=\wtil F_{i_1}^{m_1}\cd
%\wtil F_{i_l}^{m_l}(\bfk[z])
%\]
%\end{pro}
 
%%%%%%%%%%%%%%%%%%%%%%%%%%%%%%%%%%%%%%%%%%%%%%%%%%
\subsection{Localization}
\begin{df}[\cite{Loc}]
Let $M$ be a simple $R$-module. A graded braider $(M,R_M,\phi)$ 
is {\it non-degenerate} if $R_M(L(i)):M\circ L(i)\to
q^{\phi(\al_i)}L(i)\circ M$ is a non-zero homomorphism.
\end{df}
For $\rgmod$, there exists
a non-degenerate real commuting family of graded braiders
$(C_i,R_{C_i},\phi_i)_{i\in I}$(\cite{Loc}).
Set $C_\Lm:={\bf M}(w_0\Lm,\Lm)$ and denote $C_{\Lm_i}$ by 
$C_i$. 
\begin{pro}[\cite{strat}]
For $\Lm=\sum_im_i\Lm_i\in P_+$,  
we obtain the following isomorphism up to grading shift:
\begin{equation}
C_\Lm:={\bf M}(w_0\Lm,\Lm)\cong C_1^{\circ m_1}\circ \cd\circ C_n^{\circ m_n}.
\label{prod-1-n}
\end{equation}
\end{pro}
%%%%%%%%%%%%%%%%%%%%%%%%%%%%%%%%%%
\begin{thm}[{\cite[Proposition 5.1]{Loc}}]
%Let $w_0$ be the longest Weyl group element and 
Define the function $\phi_{i}:Q\to\bbZ$ by 
\[
\phi_{i}(\beta):=-(\beta,w_0\Lm_i+\Lm_i).
\]
Then there exists $\{(C_i,R_{C_i},\phi_{i})\}_{i\in I}$
 a non-degenerate real commuting family
of graded braiders of the monoidal category $\rgmod$.
\end{thm}
Now, we take $\Gamma=P=\bigoplus_i\bbZ\Lm_i$ and 
$\Gamma_+=P_+=\bigoplus_i\bbZ_{\geq0}\Lm_i$. Here, 
we obtain the localization $R\hbox{-gmod}[C_i^{\circ-1}\,|\,i\in I]$ by 
$\{(C_i,R_{C_i},\phi_{i})\}_{i\in I}$, which will be denoted by
$\tilrgmod$.

By the above Proposition, it holds the following properties:
\begin{pro}[\cite{Loc}]
Let $\Phi:\rgmod\to\tilrgmod$ be the canonical functor. Then,
\begin{enumerate}
\item
$\tilrgmod$ is an abelian category and the functor $\Phi$ is exact.
\item
For any simple object $S\in\rgmod$, $\Phi(S)$ is simple in $\tilrgmod$.
\item
$\wtil C_i:=\Phi(C_i)$ ($i\in I$) 
is invertible central graded braider in $\tilrgmod$.

\medskip

\nd
For $\mu\in P$, define $\wtil C_\mu$ such that  
$\wtil C_\mu:=\Phi(C_\mu)$ for $\mu\in P_+$, 
$\wtil C_{-\Lm_i}=C_i^{\circ -1}$ and 
$\wtil C_{\lm+\mu}=\wtil C_{\lm}\circ \wtil C_{\mu}$ for $\lm,\mu\in P$
up to grading shift.
\item
Any simple object in $\tilrgmod$ is isomorphic to $\wtil C_\Lm\circ \Phi(S)$ 
for some simple module $S\in\rgmod$ and $\Lm\in P$.
%\item
%For simple objects $S,\,S'\in\rgmod$ and $\Lm,\Lm'\in P$, 
%\[
%\wtil C_\Lm\circ\Phi(S)\cong \wtil C_{\Lm'}\circ\Phi(S')
%\Longleftrightarrow \exists \mu\in P\hbox{ s.t. }\Lm+\mu,\Lm'+\mu\in P_+
%\hbox{ and }
%q^{H(\Lm,\mu)}C_{\Lm+\mu}\circ S\cong q^{H(\Lm',\mu)}C_{\Lm'+\mu}\circ S'
%\]
\end{enumerate}
Note that in (4) $\Lm\in P$ and $S\in\rgmod$ are not necessarily unique.
\end{pro}
\begin{rem}
In \cite{Loc}, the localization is
applied to more general category $\scC_w$, which is the full subcategory
of $R$-gmod associated with a Weyl group element $w$. 
The category $R$-gmod here coincides with $\scC_{w_0}$ associated with the 
longest element $w_0$ in $W$.
\end{rem}
\begin{df}
The category $\tilrgmod$ is a $\Gamma$-graded $\bf k$-linear abelian and monoidal category. Therefore, 
its Grothendieck ring $\cK(\tilrgmod)$ holds a natural 
$\bbZ[q,q^{-1}]$-algebra structure, which defines a
{\it localized quantum coordinate ring}
$\wtil\aq:=\bbQ(q)\ot_{\bbZ[q,q^{-1}]}\cK(\tilrgmod)$.
\end{df}
Indeed, the Grothendieck ring $\cK(\tilrgmod)$ is described as follows:
\begin{df}
For a ring $R$ (not necessarily commutative) and a multiplicative set $S\subset R$, 
a ring $R'$ is said to be a {\it left ring of quotients} of $R$ with respect to $S$ if 
there exists a homomorphism $\vp:R\to R'$ such that
\begin{enumerate}
\item
Any $s'\in \vp(S)$ is invertible in $R'$.
\item
Any $m\in R'$ is in the form $m=\vp(s)^{-1}\vp(a)$ for some $s\in S,\,\,a\in R$.
\item
${\rm Ker}\,\vp=\{r\in R\mid sr=0\hbox{ for some }s\in S\}$.
\end{enumerate}
$R'$ is denoted by $S^{-1}R$.
\end{df}
\begin{pro}[{\cite[Corollary 5.4]{Loc}}]
The Grothendieck ring $\cK(\tilrgmod)$ is isomorphic to the left ring of quotients of the ring 
$\cK(\rgmod)$ with respect to the multiplicative set 
\[
\cS:=\{q^k\prod_{i\in I}[C_i]^{a_i}\mid k\in\bbZ,\,\, (a_i)_{i\in I}\in \bbZ^{I}_{\geq0}\},
\]
that is, $\cK(\tilrgmod)\cong \cS^{-1}\cK(\rgmod)$.
\end{pro}
%%%%%%%%%%%%  Sect. 7 %%%%%%%%%%%%%%%%%%%%%%%%%%%%%%%
\section{Crystal Structure on localized quantum coordinate rings
%$\wtil R$-$\rm gmod$
}
%First, we introduce:
%\begin{lem}[\cite{K-loc}]%Lemma 5.3
%For any simple module $M\in\tilrgmod$, there exists a unique $n\in\bbZ$
%such that $q^n M$ is self-dual.
%\end{lem}
In this section, one of the main theorems will be introduced, 
which claims that we can define a certain crystal structure on the localized quantum coordinate ring $\wtil\aq$, more precisely, 
a crystal structure will be defined on the family of all 
self-dual simple modules in $\tilrgmod$.

The following lemma is essential for this section.
\begin{lem}[{\cite[Proposition 2.18]{K-L}}]\label{exact-lem}
%matome-note p62
For any $i\in I$, $\beta,\gamma\in Q_+$, any modules $M\in R(\beta)$-gmod and $N\in R(\gamma)$-gmod, 
one has the following exact sequence in $R(\beta+\gamma-\al_i)$-gmod:
\begin{equation}
0\longrightarrow E_i M\circ N\longrightarrow E_i(M\circ N)\longrightarrow
q^{-(\al,\,\beta)}M\circ E_i N\longrightarrow 0.
\label{exact}
\end{equation}
\end{lem} 
For $i\in I$, let $ i^*\in I$ be a unique 
index satisfying $ \Lm_{i^*}=-w_0\Lm_i$.
\begin{lem}\label{EiCS}
\begin{enumerate}
\item
 For $S\in \rgmod$ and $i\in I$, if $E_iS=0$, then the module
$E_iC_{\Lm_{i^*}}\circ S$ is a simple module.
\item %matome note pp101
If $E_iS=0$ for $S\in \rgmod$, then we get for $\Lm\in P_+$ with $\lan h_{i^*},\Lm\ran>0$, 
\begin{equation}
{\rm soc}(E_i(C_\Lm\circ S))\cong C_{\Lm-\Lm_{i^*}}\circ(E_iC_{i^*}\circ S),
\label{soc-i*}
\end{equation}
up to grading shift.
\end{enumerate}
\end{lem}
{\sl Proof.} (1)
By applying \eqref{exact} in Lemma \ref{exact-lem} to 
$C_{\Lm_{i^*}}\circ S\cong S\circ C_{\Lm_{i^*}}$ up to grading shift, we know that
\[
E_iC_{\Lm_{i^*}}\circ S\cong E_i(C_{\Lm_{i^*}}\circ S)\cong 
E_i(S\circ C_{\Lm_{i^*}})\cong S\circ E_iC_{\Lm_{i^*}}
\]
up to grading shift.
Now it follows from  Proposition 3.13 in \cite{Loc} that
$E_iC_{\Lm_{i^*}}\circ S$ holds a simple socle 
 $S\Delta E_iC_{\Lm_{i^*}}$ and a simple head $E_iC_{\Lm_{i^*}}\nabla S$ and 
\[
S\Delta E_iC_{\Lm_{i^*}}\cong E_iC_{\Lm_{i^*}}\nabla S
\cong {\rm Im}\,\,{\bf  r}_{E_iC_{\Lm_{i^*}},\,S}\ne0,
\]
where ${\bf r}_{M,N}:M\circ N\to N\circ M$ is a R-matrix 
for $M,N\in \rgmod$ introduced in \cite[3.2]{Loc}.
By Proposition 3.13 in \cite{Loc}  we also have that 
\[
\HOM_R(E_iC_{\Lm_{i^*}}\circ S, E_iC_{\Lm_{i^*}}\circ S)=
{\bf  k}\,{\rm id}_{E_iC_{\Lm_{i^*}}\circ S},\qq
\HOM_R(E_iC_{\Lm_{i^*}}\circ S, S\circ E_iC_{\Lm_{i^*}})=
{\bf  k}\,{\bf  r}_{E_iC_{\Lm_{i^*}},\,S}.
\]
Here let us denote the isomorphism 
$\xi:S\circ E_iC_{\Lm_{i^*}}
\mapright{\sim} E_iC_{\Lm_{i^*}}\circ S$. Since 
$\xi\circ {\bf  r}_{E_iC_{\Lm_{i^*}},\,S}\in 
\HOM(E_iC_{\Lm_{i^*}}\circ S, E_iC_{\Lm_{i^*}}\circ S)$ and 
${\bf  r}_{E_iC_{\Lm_{i^*}},\,S}$ never vanishes, 
 we obtain $\xi\circ {\bf  r}_{E_iC_{\Lm_{i^*}},\,S}
=c\,{\rm id}_{E_iC_{\Lm_{i^*}}\circ S}$ for some $c\in{\bf k}^\times$, which implies that 
\[
{\rm Im}\,\,{\bf  r}_{E_iC_{\Lm_{i^*}}, S}\cong 
S\circ E_iC_{\Lm_{i^*}}\cong E_iC_{\Lm_{i^*}}\circ S,
\]
and then $E_iC_{\Lm_{i^*}}\circ S$ is a simple module in $\rgmod$.

\medskip
\nd
(2) 
Since $E_iS=0$, by \eqref{exact} one has $E_i C_\Lm\circ S\cong E_i(C_\Lm\circ S)$ and then 
\begin{equation}
{\rm soc}( E_i C_\Lm\circ S)\cong{\rm soc}(E_i(C_\Lm\circ S)).
\label{soc-soc}
\end{equation}
Let us show 
\begin{equation}
{\rm soc}(E_iC_\Lm)\cong C_{\Lm-\Lm_{i^*}}\circ E_iC_{i^*},
\label{EiEi*}
\end{equation}
up to grading shift.
By the result of (1) and \eqref{exact}, one knows that $C_{\Lm-\Lm_{i^*}}\circ E_iC_{i^*}\cong 
E_iC_{i^*}\circ C_{\Lm-\Lm_{i^*}}$ is a simple submodule 
of $E_i(C_\Lm)\cong E_i(C_{i^*}\circ C_{\Lm-\Lm_{i^*}})$. It is also known that 
$E_i(C_\Lm)$ holds a simple socle \cite{simple} and then 
we obtain \eqref{EiEi*}.
Due to (1) and \eqref{EiEi*} one finds that $C_{\Lm-\Lm_{i^*}}\circ E_iC_{\Lm_{i^*}}\circ S$ is a simple submodule
of $E_iC_\Lm\circ S$ and since $E_iC_\Lm\circ S$ possesses a simple socle, we get
\eqref{soc-i*}.     \qed

\medskip
We set 
\[
\bbB(\tilrgmod):=\{L\,|\,L\hbox{ is a self-dual simple module in }\tilrgmod\}.
\]
\begin{lem}[\cite{Loc}]\label{qn}
For any simple $L\in\tilrgmod$, there exists a unique 
$n\in\bbZ$ such that 
$q^nL$ is self-dual simple.
For a simple module $L\in\tilrgmod$ we define $\del(L)$ to be this 
integer $n$. 
\end{lem}
Then by this lemma, we find that $\bbB(\tilrgmod)$ includes all simple modules
in $\tilrgmod$ up to grading shift.
%As above for any simple 
%module is expressed  as $\wtil C_\Lm\circ \Phi(S)$ 
%For $\Lm\in P$ and a simple $S\in \bbB(\rgmod)$, where we write simply
For a simple object 
$\wtil C_\Lm\circ \Phi(S)\in\tilrgmod$ we write simply $C_\Lm\circ S$ if there is no 
confusion.

Now let us define the Kashiwara operators 
$\wtil F_i$ and $\wtil E_i$ ($i\in I$) on 
$\bbB(\tilrgmod)$ by 
\begin{eqnarray}
&&\wtil F_i(C_\Lm\circ S)=q^{\del(C_\Lm\circ \wtil F_iS)}C_\Lm\circ \wtil F_i S,\label{def-Fi}\\
&&\wtil E_i(C_\Lm\circ S)=
\begin{cases}
q^{\del(C_\Lm\circ \wtil E_iS)}C_\Lm\circ \wtil E_i S &\hbox{ if }E_i S\ne0,\\
q^{\del(C_{\Lm-\Lm_{i^*}}\circ (\wtil E_i C_{\Lm_{i^*}}\circ S))}C_{\Lm-\Lm_{i^*}}\circ (\wtil E_i C_{\Lm_{i^*}}\circ S) &\hbox{ if }E_i S=0,
\end{cases}
\label{def-Ei}
\end{eqnarray}
where $C_\Lm\circ S$ is a self-dual simple module in $\tilrgmod$,
the actions 
$\Eit S$ and $\Fit S$ are given in \eqref{eitfit}, which is defined on 
the family of all self-dual simple modules in $\rgmod$ and 
in \eqref{def-Ei} the module $\wtil E_i C_{\Lm_{i^*}}\circ S$ is simple
by Lemma \ref{EiCS}. 
Note that for any $m>0$, $\wtil E_i^m (C_\Lm\circ S)\ne0,\,\,
\wtil F_i^m(C_\Lm\circ S)\ne0$.
 
% {Well-definedness}
%For example, the well-definedness of $\wtil E_i$ is shown by using the 
%exact sequence:%Kashiwara tronto p23
%\[
%0\longrightarrow E_i(M)\circ N\longrightarrow E_i(M\circ N)\longrightarrow 
%q^{-(\al_i,\,\beta)}M\circ E_i(N)\longrightarrow 0\q\hbox{(exact)},
%\]
%where $M\in R(\beta)$-gmod, $N\in R(\gamma)$-gmod.
%Indeed, applying this to 
%$M\circ N=C_{\Lm+\mu}\circ S\cong C_{\Lm'+\mu}\circ S'$ 
%with $\Lm+\mu, \Lm'+\mu\in P_+$ and taking the head, we get well-definedness,
%{\it i.e.,}
%\[
%C_{\Lm+\mu-\Lm_{i^*}}\circ (\wtil E_iC_{\Lm_i*}\nabla S)
%\cong C_{\Lm'+\mu-\Lm_{i^*}}\circ (\wtil E_iC_{\Lm_i*}\nabla S')
%\]
%The well-definedness of $\wtil F_i$ is proved by the fact that
%$C_\mu$ is a (quasi-)center and 
%\[
%
%C_\mu\circ (L(i)\bigtriangleup S)\cong L(i)\bigtriangleup
%(C_\mu\circ S) \q(\mu\gg0)
%\]
% 
Let $\Psi:\bbB(\rgmod)\mapright{\sim}B(\ify)$ be as in Theorem \ref{LV}.
For $C_\Lm\circ S\in \bbB(\tilrgmod)$, we also define 
\begin{equation}
\begin{array}{ll}
&\vep_i(C_\Lm\circ S)=\vep_i(\Psi(S))-\lan h_i,w_0\Lm\ran,\qq
\wt(C_\Lm\circ S)=\wt(\Psi(S))+w_0\Lm-\Lm,\\
&\vp_i(C_\Lm\circ S)=\vep_i(\Psi(C_\Lm\circ S))
+\lan h_i,\wt(C_\Lm\circ S)\ran.
\end{array}
\label{ep-wt}
\end{equation}
 
\begin{thm}\label{main1}
The 6-tuple 
$(\bbB(\tilrgmod),\wt,\{\vep_i\},\{\vp_i\},
\{\wtil E_i\},\{\wtil F_i\})_{i\in I}$ 
is a crystal.
\end{thm}
{\sl Proof.}
First, we should show the
well-definedness of all data $(\wt,\{\vep_i\},\{\vp_i\},
\{\wtil E_i\},\{\wtil F_i\})_{i\in I}$, i.e.,
these data do not depend on the presentation of a simple module $L=C_\Lm\circ S\cong C_{\Lm'}\circ S'$. 
Let us assume that for a simple module $L\in \bbB(\tilrgmod)$ there exist
$\Lm,\Lm'\in P$ and $S,\,S'\in \bbB(\rgmod)$ such that 
$L=C_\Lm\circ S\cong C_{\Lm'}\circ S'$. 
By the definition of localization, one finds that 
there is a weight $\mu\in P$ such that $\Lm+\mu,\,\Lm'+\mu\in P_+$ and 
$C_{\Lm+\mu}\circ S\cong C_{\Lm'+\mu}\circ S'$ in $R$-gmod up to grading shift.
Let us show the following lemma:
\begin{lem}\label{lem-Fi}%matome-note 69
For a dominant weight $\mu\gg0$, $i\in I$ and a simple module $S\in \rgmod$, 
one has %
the following up to grading shift:
\begin{equation}
C_\mu\circ(L(i)\nabla S)\cong L(i)\nabla(C_\mu\circ S).
\label{CiS}
\end{equation}
\end{lem}
{\sl Proof.}
Since $C_\mu$ is a central object by Theorem 5.2 in \cite{Loc}, one has
\begin{equation}
C_\mu\circ L(i)\circ S\cong L(i)\circ C_\mu\circ S,
\label{CLS}
\end{equation}
up to grading shift. 
Both $C_\mu\circ(L(i)\nabla S)$ and $L(i)\nabla(C_\mu\circ S)$ are simple
quotients in $C_\mu\circ L(i)\circ S\cong L(i)\circ C_\mu\circ S$.
Since $C_\mu$ is a real simple module, from \cite{jams} and \cite{simple} we have 
that $C_\mu\circ S$ is simple and 
$L(i)\circ C_\mu\circ S$ holds a simple head and then 
it is isomorphic to $C_\mu\circ(L(i)\nabla S)$ and $L(i)\nabla(C_\mu\circ S)$, 
which implies that 
$C_\mu\circ(L(i)\nabla S)\cong L(i)\nabla(C_\mu\circ S)$
up to grading shift.\qed

By this lemma, 
we obtain 
\[
C_{\Lm+\mu}\circ \Fit S\cong C_{\Lm+\mu}\circ (L(i)\nabla S)
\cong L(i)\nabla(C_{\Lm+\mu}\circ S)\cong
L(i)\nabla(C_{\Lm'+\mu}\circ S')\cong C_{\Lm'+\mu}\circ (L(i)\nabla S')
\cong C_{\Lm'+\mu}\circ (\Fit S')
\]
up to grading shift. It shows the well-definedness of $\Fit$.

Next, let us see the one for $\Eit$. By  Lemma \ref{exact-lem}, it yields the following exact sequence:
\begin{equation}
0\longrightarrow E_i C_{\Lm+\mu}\circ S\longrightarrow
E_i (C_{\Lm+\mu}\circ S)\longrightarrow q^{(\al_i,\xi)}C_{\Lm+\mu}\circ E_i S
\longrightarrow 0,\label{exact-CS}
\end{equation}
where we set $\xi=w_0(\Lm+\mu)-(\Lm+\mu)$.
Since $C_{\Lm+\mu}\circ S$ is simple, 
by Lemma \ref{simple} one knows that 
$E_i(C_{\Lm+\mu}\circ S)$ possesses a simple head and then a unique maximal 
submodule $M$, which includes the module $E_i C_{\Lm+\mu}\circ S$.
Set $M':=M/E_i C_{\Lm+\mu}\circ S$. 
Now, we get the following commutative diagram:
\begin{equation}
\begin{array}{ccccccccc}
&&0&&0&&0&&\\
&&\downarrow&&\downarrow&&\downarrow&&\\
0&\to&E_i C_{\Lm+\mu}\circ S&\mapright{=} &E_i C_{\Lm+\mu}\circ S&\to&0&\to&0\\
&&\downarrow&&\downarrow&&\downarrow&&\\
0&\to&M&\mapright{\,h\,} &E_i (C_{\Lm+\mu}\circ S)&\to&{\rm hd}(E_i (C_{\Lm+\mu}\circ S))&\to&0\\
&&\downarrow&&\q\downarrow\psi&&\downarrow&&\\
0&\to&M'&\mapright{\,g\,} & q^{(\al_i,\xi)}C_{\Lm+\mu}\circ E_i S&\to&
q^{(\al_i,\xi)}(C_{\Lm+\mu}\circ E_i S)/M'&\to&0\\
&&\downarrow&&\downarrow&&\downarrow&&\\
&&0&&0&&0&&
\end{array}.
\end{equation}

In this diagram, all of three rows are exact and the first two columns are also exact.
Therefore, applying the nine lemma to this diagram, we obtain 
\begin{equation}
{\rm hd}(E_i (C_{\Lm+\mu}\circ S))\cong q^{(\al_i,\xi)}
(C_{\Lm+\mu}\circ E_i S)/M'.\label{col}
\end{equation}
\begin{lem}
The module $M'$ is a unique maximal submodule of 
$q^{(\al_i,\xi)}C_{\Lm+\mu}\circ E_i S$.
% See Note22 p96-97
\end{lem}
{\sl Proof.}
The maximality of $M'$ is trivial from \eqref{col}.
%First let us show the maximality of $M'$.
%Assume that there exists a submodule $\ovl M$ such that
%$M'\cong g(M')\subsetneqq \ovl M\subsetneqq q^{(\al_i,\xi)}C_{\Lm+\mu}\circ E_iS$.
%Since $\psi^{-1}(\ovl M)$ is a submodule of $E_i(C_{\Lm+\mu}\circ S)$ and 
%$\psi\circ h(M)\cong g(M')$, we have $h(M)\subsetneqq \psi^{-1}(\ovl M)$.
%But, since $h(M)$ is a maximal submodule of $E_i(C_{\Lm+\mu}\circ S)$, one has  
%$\psi^{-1}(\ovl M)=E_i(C_{\Lm+\mu}\circ S)$ and then 
%$M'=q^{(\al_i,\xi)}C_{\Lm+\mu}\circ E_iS$, which derives a contradiction. 
%\medskip
Let us see the uniqueness of $M'$.
Suppose that there exists a maximal submodule $\wtil M\subsetneqq   
q^{(\al_i,\xi)}C_{\Lm+\mu}\circ E_iS$. It is evident that 
$\psi^{-1}(\wtil M)$ is a maximal submodule of $E_i(C_{\Lm+\mu}\circ S)$.
Since $h(M)$ is a unique maximal submodule of $E_i(C_{\Lm+\mu}\circ S)$,
we get $\psi^{-1}(\wtil M)=h(M)$ and then 
$\wtil M=\psi\circ h(M)=g(M')\cong M'$.
\qed

\medskip
%%%%%%%%%%%%%%%%%%%%%%%
By this lemma one has that $M'$ is a unique maximal submodule of 
$q^{(\al_i,\xi)}C_{\Lm+\mu}\circ E_i S$
% See Note22 p96-97
and then 
${\rm hd}(q^{(\al_i,\xi)}C_{\Lm+\mu}\circ E_i S) 
\cong q^{(\al_i,\xi)}(C_{\Lm+\mu}\circ E_i S)/M'$. 
Therefore, one obtains 
\[
{\rm hd}(E_i (C_{\Lm+\mu}\circ S))
\cong {\rm hd}(q^{(\al_i,\xi)}C_{\Lm+\mu}\circ E_i S).
\]
Hence, for $C_{\Lm+\mu}\circ S\cong C_{\Lm'+\mu}\circ S'$ 
(up to grading shift)
it follows that 
\[
{\rm hd}(C_{\Lm+\mu}\circ E_iS)\cong 
{\rm hd}(E_i(C_{\Lm+\mu}\circ S))\cong 
{\rm hd}(E_i(C_{\Lm'+\mu}\circ S'))\cong 
{\rm hd}(C_{\Lm'+\mu}\circ E_iS'),
\]
up to grading shift.
We also have 
\begin{equation}
C_{\Lm+\mu}\circ \Eit S\cong{\rm hd}(C_{\Lm+\mu}\circ E_iS),\qq
C_{\Lm'+\mu}\circ \Eit S'\cong{\rm hd}(C_{\Lm'+\mu}\circ E_iS'),
\label{Lm-Lm'}
\end{equation}
up to grading shift.
Thus, if $E_i S,E_i S'\ne0$, we obtain 
$C_{\Lm+\mu}\circ \Eit S\cong C_{\Lm'+\mu}\circ \Eit S'$
up to grading shift.
Let us assume that $E_i S=0$. 
Owing to \eqref{exact} one has 
$E_i C_{\Lm+\mu}\circ S\cong E_i(C_{\Lm+\mu}\circ S)$. By \cite[Proposition 3.23]{Loc}
one has 
\begin{equation}
\vep_i(C_j)=\begin{cases}1&\hbox{ if }j=i^*,\\
0&\hbox{ otherwise.}
\end{cases}
\label{vepi-C}
\end{equation}
Since 
$E_i C_{\Lm_{i^*}}\cong \wtil E_i C_{\Lm_{i^*}}$ by \cite[Corollary 10.1.7]{jams},
one has
\begin{equation}
{\rm  soc}(E_i C_{\Lm+\mu})\cong C_{\Lm+\mu-\Lm_{i^*}}\circ E_i C_{\Lm_{i^*}}
\cong C_{\Lm+\mu-\Lm_{i^*}}\circ \Eit C_{\Lm_{i^*}},
\label{EC}
\end{equation}
up to grading shift.
Hence, if $E_i S'=0$, one gets 
\begin{eqnarray*}
C_{\Lm+\mu-\Lm_{i^*}}\circ \Eit C_{\Lm_{i^*}}\circ S&\cong &
C_{\Lm+\mu-\Lm_{i^*}}\circ E_i C_{\Lm_{i^*}}\circ S\cong 
{\rm soc}(E_iC_{\Lm+\mu})\circ S\cong {\rm soc}(E_i(C_{\Lm+\mu}\circ S))\\
&\cong &{\rm soc}(E_i(C_{\Lm'+\mu}\circ S'))\cong
{\rm soc}(E_iC_{\Lm'+\mu})\circ S'\cong
C_{\Lm'+\mu-\Lm_{i^*}}\circ E_i C_{\Lm_{i^*}}\circ S'\\
&\cong&
C_{\Lm'+\mu-\Lm_{i^*}}\circ \Eit C_{\Lm_{i^*}}\circ S',
\end{eqnarray*}
up to grading shift.
Even if $E_i S'\ne0$, it follows from \eqref{eitfit}, \eqref{Lm-Lm'} and \eqref{EC}
that the simple module $C_{\Lm'+\mu}\circ \Eit S'$ is isomorphic to 
\[
{\rm soc}E_i(C_{\Lm'+\mu}\circ S')\cong {\rm soc}E_i(C_{\Lm+\mu}\circ S)
\cong C_{\Lm+\mu-\Lm_{i^*}}\circ (\Eit C_{\Lm_{i^*}}\circ S)
\]
up to grading shift, which implies
the well-definedness of $\Eit$ on $\bbB(\tilrgmod)$.

Let us investigate the well-definedness of the function $\vep_i$.
For $n>0$, $\mu\in P_+$ and a simple module $S\in R$-gmod, 
the following formula is an immediate consequence of \eqref{exact-CS}:
\begin{equation}%see matome-note p73
[E_i^n(C_\mu\circ S)]=\bigoplus_{k=0}^na_k(q)[E_i^k C_\mu\circ E_i^{n-k}S],
\label{sum-Ei}
\end{equation}
where $[M]$ means an equivalence class of a module $M$ in Grothendieck ring ${\cK}(\rgmod)$ and 
$a_k(q)$ is a non-zero 
Laurent polynomial in $q$ with non-negative coefficients. 
It implies that $\vep_i(C_\mu\circ S)=\vep_i(C_\mu)+\vep_i(S)$.
Therefore, since $\vep_i(C_{\Lm_j^*})=\del_{i,j}$ and then $\vep_i(C_\mu)$ coincides with the number of $C_{\Lm_{i^*}}
=C_{-w_0\Lm_i}$ in $C_\mu$, one gets $\vep_i(C_\mu)=\lan h_i,-w_0\mu\ran$. Thus, 
\[\vep_i(C_\mu\circ S)=\vep_i(S)-\lan h_i,w_0\mu\ran.
\]
Hence, if $C_{\Lm+\mu}\circ S\cong C_{\Lm'+\mu}\circ S'$ 
for $\Lm+\mu,\Lm'+\mu\in P_+$, then
$\vep_i(C_{\Lm+\mu}\circ S)=\vep(C_{\Lm+'\mu}\circ S')$ and then 
\[
\vep_i(S)-\lan h_i,w_0(\Lm+\mu)\ran=\vep_i(S')-\lan h_i,w_0(\Lm'+\mu)\ran,
\]
holds, which shows 
$\vep_i(S)-\lan h_i,w_0\Lm\ran=\vep_i(S')-\lan h_i,w_0 \Lm'\ran$.  
Now, we get the well-definedness of $\vep_i$.
The well-definedness of $\vp_i$ and $\wt$ is trivial, which 
completes the proof for the well-definedness of $\Eit,\Fit,\vep_i,\vp_i$ and $\wt$
on $\bbB(\tilrgmod)$.

Next, let us see that the 6-tuple satisfies the conditions of crystal
as in Definition \ref{cryst} (1)--(5). The condition (1) is automatically satisfied
by the definition of $\vp_i$. As for the condition (5), $\vp_i(b)=-\ify$ never happen.
So we shall see the condition (2). 
For $C_\Lm\circ S\in\bbB(\tilrgmod)$, suppose that $\Eit S\ne0$.
\[
\vep_i(\Eit(C_\Lm\circ S))=\vep_i(C_\Lm\circ \Eit S)
=\vep_i(\Eit S)-\lan h_i,w_0\Lm\ran
=\vep_i(S)-1-\lan h_i,w_0\Lm\ran=\vep_i(C_\Lm\circ S)-1.
\]
It is trivial that $\wt(C_\Lm\circ \Eit S)=\wt(\Eit S)+w_0\Lm-\Lm= %shuusei 
\wt(S)+\al_i+w_0\Lm-\Lm=\wt(C_\Lm\circ S)+\al_i$.   %shuusei 
Suppose that $E_i S=0$, which implies that 
$0=\vep_i(S)=\vep_i(\wtil E_iC_{\Lm_{i^*}}\circ S)$ 
and $\vep_i(C_\Lm\circ S)=
\vep_i(C_\Lm)+\vep_i(S)=-\lan h_i,w_0\Lm\ran$. Then, we have
\[
\vep(\Eit(C_\Lm\circ S))=\vep_i(C_{\Lm-\Lm_{i^*}}\circ  (\Eit C_{\Lm_{i^*}}\circ S))
=\vep_i(\Eit C_{\Lm_{i^*}}\circ S)-\lan h_i,w_0({\Lm-\Lm_{i^*}})\ran
=-\lan h_i,w_0({\Lm})\ran -1
=\vep_i(C_\Lm\circ S)-1.
\]%matome note p19
We also get 
\begin{eqnarray*}
\wt(C_{\Lm-\Lm_{i^*}}\circ  (\Eit C_{\Lm_{i^*}}\circ S))&=&
\wt(\Eit C_{\Lm_{i^*}}\circ S)+w_0(\Lm-\Lm_{i^*})-\Lm+\Lm_{i^*}\\
%\qq=\wt(\Eit C_{\Lm_{i^*}}\circ S)+w_0(\Lm-\Lm_{i^*})-\Lm+\Lm_{i^*}
&=&\wt(S)+\wt(\Eit C_{\Lm_{i^*}})+w_0(\Lm-\Lm_{i^*})-\Lm+\Lm_{i^*}\\
&=&\wt(S)+w_0(\Lm_{i^*})-\Lm_{i^*}+\al_i+w_0(\Lm-\Lm_{i^*})-\Lm+\Lm_{i^*}\\
&=&\wt(S)+w_0\Lm-\Lm+\al_i=\wt(C_\Lm\circ S)+\al_i
\end{eqnarray*}

The condition (3) will be shown by similar way to the condition (2).
Let us check the condition (4). To do that, it suffices to show that 
$\Eit\Fit={\rm id}_{\bbB(\wtil R{\rm -gmod})}=\Fit\Eit$.
At first, since $\Eit(\Fit S)=S\ne0$ for any $S\in \bbB(\rgmod)$,
it is trivial that $\Eit\Fit(C_\Lm\circ S)=
\Eit(C_\Lm\circ \Fit S)=C_\Lm\circ S$. 
In the case $E_i S\ne0$, it is also evident that
$\Fit\Eit(C_\Lm\circ S)=\Fit(C_\Lm\circ \Eit S)\cong C_\Lm\circ \Fit\Eit S
\cong C_\Lm\circ S$.
So, assuming $E_i S=0$, let us show the following lemma:
\begin{lem}
For any simple $S\in \rgmod$ and $i\in I$, if $ E_i S=0$, then we obtain
\begin{equation}
\Fit(\Eit C_{\Lm_{i^*}}\circ S)\cong C_{\Lm_{i^*}}\circ S,
\label{FE}
\end{equation}
up to grading shift.
\end{lem}
{\sl Proof.}
It follows from \eqref{exact} that there exists the following 
exact sequence
\[
0\longrightarrow (E_i C_{\Lm_{i^*}})\circ S\longrightarrow 
E_i(C_{\Lm_{i^*}}\circ S)\longrightarrow
q^{-(\al_i,\Lm_{i^*}-w_0\Lm_{i^*})}C_{\Lm_{i^*}}\circ E_iS
\longrightarrow 0.
\]
Since $E_i S=0$ and $E_i C_{\Lm_{i^*}}\cong \Eit C_{\Lm_{i^*}}$, one has 
$(\Eit C_{\Lm_{i^*}})\circ S\cong E_i(C_{\Lm_{i^*}}\circ S)$.
Since $(\Eit C_{\Lm_{i^*}})\circ S$ is simple by Lemma \ref{EiCS}, 
one obtains 
$\Eit(C_{\Lm_{i^*}})\circ S\cong \Eit(C_{\Lm_{i^*}}\circ S)$ in $\rgmod$ 
up to grading shift, 
which shows \eqref{FE}.\qed

Therefore, by this lemma under the assumption $E_i S=0$
we get 
\[
\Fit\Eit(C_\Lm\circ S)\cong
\Fit(C_{\Lm-\Lm_{i^*}}\circ (\Eit C_{\Lm_{i^*}}\circ S))
\cong C_{\Lm-\Lm_{i^*}}\circ \Fit(\Eit C_{\Lm_{i^*}}\circ S)
\cong C_{\Lm-\Lm_{i^*}}\circ C_{\Lm_{i^*}}\circ S
\cong C_\Lm\circ S.
\]
Thus, we get that the condition (4) holds and then 
we have completed the proof of the theorem.\qed

%%%%%%%%%%%%%%%%%%%%%%%%%%%%%%%%%%%%%%%%%%%
\section{Cellular Crystal $\bbB_\bfi$ and $\bbB(\wtil R${\rm-gmod}$)$\q}

In this section, we shall describe the explicit crystal structure on 
$\bbB(\tilrgmod)$ %defined -->shuusei
by constructing the isomorphism
to the cellular crystal $\bbB_\bfi$.

Here, by Proposition \ref{pro-f0} 
we observe that there seems to exist a certain correspondence:
\begin{eqnarray*}
\{C_\Lm\,|\,\Lm\in P_+\}\subset \rgmod&\longleftrightarrow&
\cH_\bfi\\
C_\Lm=\wtil F_{i_1}^{m_1}\cd
\wtil F_{i_N}^{m_N}{\bf 1}
&\longleftrightarrow&{\bf h}_\Lm=
\til f_{i_1}^{m_1}\til f_{i_2}^{m_2}\cd \til f_{i_N}^{m_N}
((0)_{i_1}\ot(0)_{i_2}\ot\cd\ot (0)_{i_N})
\end{eqnarray*}
Together with the result of Proposition \ref{pro-bi}, we obtain the following:
%Therefore, we obtain:
\begin{thm}\label{iso}
For any reduced longest word $\bfi=i_1i_2\cd i_N$, 
there exists an isomorphism of crystals:
\begin{eqnarray*}
\wtil\Psi:\bbB(\tilrgmod)&\mapright{\sim}&\bbB_\bfi=\bigcup_{h\in\cH_\bfi}
B^{h}(\ify)\\
C_\Lm\circ S&\longmapsto& \bfh_\Lm+\Psi(S)\in B^{\bfh_\Lm}(\ify),
\end{eqnarray*}
where $\Psi:\bbB(\rgmod)\mapright{\sim} B(\ify)$ is the isomorphism 
of crystals given in Theorem \ref{LV}, 
$S$ is simple in $\bbB(\rgmod)$ and 
for $\Lm=\sum_ia_i\Lm_i$ set 
$\bfh_\Lm=\sum_ia_i\bfh_i$.
\end{thm}

 {\sl Proof.}
We shall show:
\begin{enumerate}
\item
Well-definedness of $\wtil\Psi$, that is, 
if $C_\Lm\circ S\cong C_{\Lm'}\circ S'$, then 
$\bfh_\Lm+\Psi(S)=\bfh_{\Lm'}+\Psi(S')$.
\item
Bijectivity of $\wtil\Psi$.
\item 
$\wt(C_\Lm\circ S)=\wt(\bfh_\Lm+\Psi(S))$.
\item 
$\vep_i(C_\Lm\circ S)=\vep_i(\bfh_\Lm+\Psi(S))$.
\item
$\wtil\Psi\circ \Fit=\fit\,\circ \wtil\Psi$ and 
$\wtil\Psi\circ \wtil E_i=\eit\,\circ \wtil\Psi$.
\end{enumerate}

%matome-note p65
First, let us see (1). We may show that 
${\bf h}_{\Lm+\mu}+\Psi(S)={\bf h}_{\Lm'+\mu}+\Psi(S')$
in case of $C_{\Lm+\mu}\circ S\cong 
C_{\Lm'+\mu}\circ S'$ for 
$\mu\gg0$ satisfying $\Lm+\mu,\Lm'+\mu\in P_+$. Thus, 
we shall show if $S'\cong C_\Lm\circ S$ for simple modules $S,S'\in \rgmod$
and $\Lm\in P_+$, then $\Psi(S')=\Psi(S)+{\bf h}_\Lm$.
Here we denote $(0)_{i_1}\ot\cd\ot (0)_{i_N}$ by $\bf0_\bfi$ and identify 
$u_\ify$ with $\bf 0_\bfi$.
By the definition of the determinantial module
$C_\Lm={\bf M}(w_0\Lm,\Lm)$ and the isomorphism of crystals
$\Psi:\rgmod\mapright{\sim} B(\ify)$, 
we know that there exists a sequence of indices 
$j_1,\cd,j_l,j_{l+1},\cd,j_{l+m}$ such that
$\wtil E_{j_{l+m}}\cd \wtil E_{j_{l+1}}C_\Lm={\bf 1}$ and 
$\wtil E_{j_{l}}\cd \wtil E_{j_1}S={\bf 1}$. Then, we have 
\begin{eqnarray*}
\wtil E_{j_{l+m}}\cd \wtil E_{j_{1}}(C_\Lm\circ S)
&\cong& 
\wtil E_{j_{l+m}}\cd \wtil E_{j_{l+1}}(C_\Lm\circ \wtil E_{j_{l}}\cd \wtil E_{j_1}S)\\
&\cong& 
\wtil E_{j_{l+m}}\cd \wtil E_{j_{l+1}}(C_\Lm\circ {\bf1})
\cong
\wtil E_{j_{l+m}}\cd \wtil E_{j_{l+1}}(C_\Lm)
\cong {\bf1},
\end{eqnarray*}
which also means  $\wtil E_{j_{l+m}}\cd \wtil E_{j_{1}} S'={\bf 1}$.
Then, since $\Psi:\rgmod\to B(\ify)$ is an isomorphism of crystals, 
we find that $\Psi(S')=\til f_{j_1}\cd\til f_{j_{l+m}}{\bf 0}_\bfi$.
It follows from Proposition \ref{pro-bi} (1) and Proposition \ref{pro-f0}
that
\[
\til f_{j_1}\cd\til f_{j_{l+m}}{\bf 0}_\bfi=
\til f_{j_1}\cd\til f_{j_l}(\til f_{j_{l+1}}
\cd \til f_{j_{l+m}}{\bf 0}_\bfi)
=\til f_{j_1}\cd\til f_{j_l}({\bf h}_\Lm)
=\til f_{j_1}\cd\til f_{j_l}({\bf0}_\bfi+{\bf h}_\Lm)
=\Psi(S)+{\bf h}_\Lm.
\] 
Hence, we obtain $\Psi(S')=\Psi(S)+{\bf h}_\Lm$ and then (1).

%matome-note p51-
Next, let us show (2) the bijectivity of $\wtil\Psi$.
For any $x\in \bbB_\bfi$ by Proposition \ref{pro-bi}, there exists
$\Lm\in P$ and $b\in B(\ify)$ such that $x={\bf h}_\Lm+b\in B^{{\bf h}_\Lm}(\ify)$.
Then, the fact $\wtil\Psi(C_\Lm\circ \Psi^{-1}(b))={\bf h}_\Lm+b$ shows the 
surjectivity of $\wtil\Psi$.
The following lemma is almost evident by the property of the localization method:
\begin{lem} %matome-note p61
For any simple modules $M,M'\in \tilrgmod$, there exist $\xi\in P$ and simple 
modules $L,L'\in \rgmod$ such that 
$M=C_\xi\circ L$ and $M'=C_\xi\circ L'$.
\end{lem}
{\sl Proof.} 
For $M,\,M'$ there exist $\Lm,\,\Lm'\in P$ and $S,\,S'\in\rgmod$ such that
$M\cong C_\Lm\circ S$ and $M'\cong C_{\Lm'}\circ S'$.
Then, there exists a dominant weight $\mu\gg0$ such that
$\Lm+\mu,\,\Lm'+\mu\in P_+$ and 
$L:=C_\mu\circ M\cong C_{\Lm+\mu}\circ S,\,L':=C_{\Lm'+\mu}\circ S'\cong
C_\mu\circ M'\in \rgmod$ up to grading shift. Note that both $L,\,L'$ are simple.
Thus, setting $\xi=-\mu$, since 
$C_\mu$ is invertible, we get 
the desired result.\qed

For any $M,M'\in\bbB(\tilrgmod)$, by the above lemma we can set
$M=C_\mu\circ L$ and $M'=C_\mu\circ L'$ for some 
$\mu\in P$, $L,L'\in \bbB(\rgmod)$. 
Here we assume $\wtil\Psi(C_\mu\circ L)=\wtil\Psi(C_\mu\circ L')$.
Then we obtain ${\bf h}_\mu+\Psi(L)={\bf h}_\mu+\Psi(L')$, 
which implies $\Psi(L)=\Psi(L')$ and then $L\cong L'$ since $\Psi$ 
is bijective.
Then, one has $M\cong M'$, which means the injectivity of $\wtil\Psi$.

Let us see (3).
By the definition of $\wt$ in \eqref{ep-wt} and Proposition \ref{pro-f0}, 
it is immediate that 
$\wt(C_\Lm\circ S)=w_0\Lm-\Lm+\wt(\Psi(S))=\wt({\bf h}_\Lm+\Psi(S))$. 

%matome-note p22
As for (4), first by the definition of $\vep_i$ on $\bbB(\tilrgmod)$
in \eqref{ep-wt}
one has
$\vep_i(C_\Lm\circ S)=\vep_i(\Psi(S))-\lan h_i,w_0\Lm\ran$.
Now, set ${\bf  h}_\Lm=(m_1,\cd,m_N)\in \cH_\bfi$ and $\Psi(S)=(x_1,\cd,x_N)\in B(\ify)$, where 
$m_k:=\lan h_{i_k},s_{i_{k+1}}\cd s_{i_N}\Lm\ran$.
Let $\{k_1,k_2,\cd\}=\{k\mid i=i_k\}$ where $k_1<k_2<\cd$. 
By the definition of $\cH_\bfi$, one has $\sigma_{k_1}({\bf  h}_\Lm)
=\sigma_{k_2}({\bf  h}_\Lm)=\cd$ and then
\[
\vep_i({\bf  h}_\Lm)=\sigma_{k_1}({\bf h}_\Lm)=m_{k_1}+\sum_{j<k_1}a_{i,i_j}m_j.
\]
Using this formula, one gets
\begin{eqnarray*}
\vep_i({\bf  h}_\Lm+\Psi(S))
&=&\max\{\sigma_k({\bf  h}_\Lm+\Psi(S))\mid i=i_k\}
=\max\{x_k+m_k+\sum_{j<k}a_{i,i_j}(x_j+m_j)\mid i=i_k\}\\
&=&\max\{x_k+\sum_{j<k}a_{i,i_j}x_j\mid i=i_k\}+
m_{k_1}+\sum_{j<k_1}a_{i,i_j}m_j
=\vep_i(\Psi(S))+\vep_i({\bf  h}_\Lm).
\end{eqnarray*}
It yields that %matome-note p23
\begin{eqnarray*}
\lan h_i,w_0\Lm\ran&=&\lan h_i, s_{i_1}\cd s_{i_N}\Lm\ran
=\lan h_i,s_{i_2}\cd s_{i_N}\Lm-\lan h_{i_1},s_{i_2}\cd s_{i_N}\Lm\ran\al_{i_1}\ran
\\
& = &
\lan h_i,s_{i_2}\cd s_{i_N}\Lm\ran-\lan h_{i_1},\al_{i_2}\ran m_1
=\lan h_i,s_{i_3}\cd s_{i_N}\Lm-\lan h_{i_2},s_{i_2}\cd s_{i_N}\Lm\ran \al_{i_2}\ran
-\lan h_i,\al_{i_1}\ran m_1
\\
&=&
\lan h_i,s_{i_3}\cd s_{i_N}\Lm\ran-\lan h_i,\al_{i_2}\ran m_2
-\lan h_i,\al_{i_1}\ran m_1
=\cd
=\lan h_i,s_{i_{k_1+1}}\cd s_{i_N}\Lm\ran-
\sum_{j\leq k_1}\lan h_i,\al_{i_j}\ran m_j\\
&=& m_{k_1}-2m_{k_1}-\sum_{j< k_1}\lan h_i,\al_{i_j}\ran m_j
=-m_{k_1}-\sum_{j< k_1}\lan h_i,\al_{i_j}\ran m_j=-\vep_i({\bf h}_\Lm),
\end{eqnarray*}
which shows (4), that is,
\[
\vep_i(C_\Lm\circ S)=\vep_i(\Psi(S))-\lan h_i,w_0\Lm\ran
=\vep_i(\Psi(S))+\vep_i({\bf h}_\Lm)=\vep_i({\bf h}_\Lm+\Psi(S)).
\]

Let us investigate (5). %matome-note p20
For a simple $C_\Lm\circ S\in \bbB(\tilrgmod)$, it yields
\[
\wtil\Psi\Fit(C_\Lm\circ S)=\wtil\Psi(C_\Lm\circ \Fit S)=
{\bf h}_\Lm+\Psi(\Fit S)={\bf h}_\Lm+\fit\Psi(S),\]
and by \eqref{eifih} in Proposition \ref{pro-bi}, one has
\[
\fit\wtil\Psi(C_\Lm\circ S)=\fit({\bf h}_\Lm+\Psi(S))
={\bf h}_\Lm+\fit\Psi(S).\]
It shows $\wtil\Psi\circ \Fit=\fit\,\circ \wtil\Psi$.
By use of the formula $\wtil\Psi\circ \Fit=\fit\,\circ \wtil\Psi$ and 
$\Fit\Eit=\Eit\Fit={\rm id}_{\bbB(\wtil R{\rm -gmod})}$,  
it is clear that
\[
{\bf h}_\Lm+\Psi(S)=\wtil\Psi(\Fit\Eit(C_\Lm\circ S))
=\fit\,\wtil\Psi(\Eit(C_\Lm\circ S))\]
and then applying $\eit$ on the both sides, it 
follows from $\fit\eit=\eit\fit={\rm id}_{\bbB_\bfi}$ that
\[
\eit({\bf h}_\Lm+\Psi(S))=\wtil\Psi(\Eit(C_\Lm\circ S)),
\]
which implies $\eit\wtil\Psi=\wtil\Psi\Eit$.
Now, we have completed the proof of the theorem.\qed

%%%%%%%%%%%%%%%%%%%%%%%%%%%%%%%%%%%%%%%%%
\section{Application and further problems}
In this section, we shall see several miscellaneous results as applications
of the crystal structure on $\bbB(\tilrgmod)$.
%\topmargin{-20pt}
%Here we list sveral further problems.
% {Problem 2}
%Relation to cluster algebras: Indeed, the determinantial modules correspond to
%frozen variable of cluster algebra.
% 

%Let us introduce the following {anti-automorphism} $*$ of $\uq$:
%\begin{df}%[{[K1],[K2]}]
\subsection{Operator $\til\fra$}
Define the $\bbQ(q)$-linear anti-automorphism $\star$ of $\uq$ by 
\begin{equation*}
(q^h)^\star=q^{-h},\q e_i^\star=e_i,\q f_i^\star=f_i.
\label{star}
\end{equation*}
%\end{df}
%Note that 
%$*^2={\rm id}$. 
\begin{thm}[\cite{K3}]\label{star-dual}
Set $L^\star(\ify):=\{u^\star\,|\,u\in L(\ify)\}$, 
$B^\star(\ify):=\{b^\star\,|\,b\in B(\ify)\}$. Then we have 
\begin{eqnarray*}
L^\star(\ify)=L(\ify),
&B^\star(\ify)=B(\ify).
\end{eqnarray*}
\end{thm}
From the proof of Theorem 5.13 in \cite{Loc} we get
\begin{pro}[\cite{Loc}]\label{a-1}
For $\nu=(\nu_1,\nu_2,\cd,\nu_{m-1},\nu_m)\in I^\beta$ ($m:=|\beta|$)
set $\ovl\nu=(\nu_m,\nu_{m-1},\cd,\nu_2,\nu_1)$. 
Define the automorphism ${\fra}$ on $R(\beta)$  by 
\[
{\fra}(e(\nu))=e(\ovl \nu),\q
{\fra}(x_ie(\nu))=x_{m-i+1}e(\ovl \nu),\q
{\fra}(\tau_je(\nu))=-\tau_{m-j}e(\ovl\nu).
\]
Then, there exists the functor $\fra:\rgmod\to\rgmod$ such that
$\fra(C_i)=C_{i^*}$ ($\forall i\in I$), $\fra^2\cong {\rm id}$ 
and  $\fra(X\circ Y)\cong \fra(Y)\circ\fra(X)$ for $X,Y\in\rgmod$. Furthermore, 
it is extended to the functor $\til\fra:\tilrgmod\to\tilrgmod$ which satisfies
\begin{equation}
\til\fra^2\cong{\rm id},\q \hbox{and}\q
\til\fra(X\circ Y)\cong \til\fra(Y)\circ\til\fra(X)\q
\hbox{ for }X,Y\in\tilrgmod.
\label{fra}
\end{equation}
\end{pro}
Note that $\fra$(resp. $\til\fra$) induces the operation $\star$ 
on $\aq$ (resp. $\tilaq$) since 
$\fra(L(i))=L(i)$ and then one has $\fra(f_i)=f_i$ (resp. $\til\fra(f_i)=f_i$) 
on $\aq$ (resp. $\tilaq$).
Now, we obtain the following:
\begin{pro}\label{a-bbB}
Let $\til\fra:\tilrgmod\to\tilrgmod$ be the functor as above.
It yields 
\begin{equation}
\til\fra(\bbB(\tilrgmod))=\bbB(\tilrgmod).
\end{equation}
\end{pro}
{\sl Proof.}
It is evident that 
for any simple module $L\in\tilrgmod$, $\til\fra(L)$ is also simple. 
Therefore, it suffices to show that 
if a simple module $L$ satisfies $L\cong L^*$, then $\til\fra(L)\cong \til\fra(L)^*$.
To do that, we show the following lemma:
\begin{lem}\label{psi-a}
Let $\psi$ be the anti-automorphism on $R(\beta)$ as in \ref{defQHA}. We obtain
\begin{equation}
\psi\circ \fra=\fra\circ \psi.
\end{equation}
\end{lem}
{\sl Proof.}  %Lemma \ref{psi-a}
For $\beta$ with $|\beta|=m$, it is trivial that for the generators $x_i, e(\nu)$ of $R(\beta)$, 
\[
\psi\circ \fra(e(\nu))=e(\ovl\nu)=\fra\circ\psi(e(\nu)),\qq
\psi\circ \fra(x_ie(\nu))=e(\ovl\nu)x_{m+1-i}=x_{m+1-i}e(\ovl\nu)
=\fra\circ\psi(x_i(e(\nu))
\]
One also gets
\begin{eqnarray*}
\psi\circ \fra(\tau_je(\nu))&=&-\psi(\tau_{m-j}e(\ovl\nu))=-e(\ovl\nu)\tau_{m-j}
=-\tau_{m-j}e(s_{m-j}(\ovl\nu))=-\tau_{m-j}e(\ovl{s_j\nu})\\
&=&\fra(\tau_je(s_j\nu))=\fra(e(\nu)\tau_j)=\fra\circ \psi(\tau_je(\nu)),
\end{eqnarray*}
which completes the proof.\qed

By this lemma, we find that 
\begin{equation}
\til\fra(L)^*\cong \til\fra(L^*).
\label{a*=*a}
\end{equation}
Since $\til\fra(L)$ is simple, by Lemma \ref{qn} there exists a unique $n\in \bbZ$ such that 
$(q^n\til\fra(L))^*\cong q^n\til\fra(L)\cong \til\fra(q^n L)$ and then 
$q^{-n}(\til\fra(L))^*\cong (q^n\til\fra(L))^*\cong\til\fra(q^n L)$. 
Therefore, owing to \eqref{a*=*a} we get 
\[
\til\fra(q^nL)\cong q^{-n}\til\fra(L^*)\cong \til\fra(q^{-n} L^*) 
\q\hbox{ and then }\q q^{2n}L\cong L^*,
\]
which means $n=0$ and then we obtain $\til\fra(L)\cong\til\fra(L)^*$.\qed

Here note that Proposition \ref{a-bbB} can be seen as a generalization of 
Theorem \ref{star-dual}.

\medskip
Since as crystals $\bbB(\tilrgmod)\cong \bbB_\bfi$ for any reduced longest word $\bfi$, 
the proposition above gives rise to the following problem.\\
\nd{\bf Problem 1.}
Can we describe $\til\fra$-operation on $\bbB_\bfi
=B_{i_1}\ot \cd\ot B_{i_N}$ explicitly?
 
 Of course, this problem is non-trivial since even for the case $B(\ify)$ the 
 explicit description has not yet been done before in $\bbB_\bfi$.
 
%\medskip

\medskip
\subsection{Summation $\oplus_\bfi$}
By use of the map $\wtil\Psi$ we can define a certain additive structure on 
$\bbB(\tilrgmod)$ as follows:
For $L,L'\in \bbB(\tilrgmod)$ and a reduced longest word $\bfi$, we define
\begin{equation}
L\oplus_\bfi L':=\wtil\Psi^{-1}(\wtil\Psi(L)+\wtil\Psi(L')),
\label{plus}
\end{equation}
where the summation in the right hand side of \eqref{plus} is the natural 
summation of $\bbB_\bfi$ by identifying with the free-$\bbZ$ lattice $\bbZ^N$.
The following is trivial from the definition of $\oplus_\bfi$
\begin{pro}
By the summation $\oplus_\bfi$, $\bbB(\tilrgmod)$ becomes an additive group.
Here let us denote the inverse of $L$ with respect to $\oplus_\bfi$ by $ L^{\ominus_\bfi}$, 
where we define $L^{\ominus_\bfi}:=\wtil\Psi^{-1}(-\wtil\Psi(L))$.
\end{pro}
Note that this summation DOES depend on a choice of reduced word $\bfi$. In this sense, 
it is not canonical. The following example implies that.

\begin{ex}%note 22 p99
For $\ge=\ssl_3$ and $\bfi=121$, let us see the inverse of $L(1)\in R(\al_1)$-gmod, which 
corresponds to $(0,0,1)\in\bbB_{121}$.
As has seen in Proposition \ref{basis-h}, 
${\bf h}_1=(0,1,1) $ and ${\bf h}_2=(1,1,0)$, and then we obtain 
\[
(0,0,-1)=\til f_2(0,-1,-1)=\til f_2(-{\bf h}_1)=-{\bf h}_1+\til f_2(0,0,0).
\]
This implies that 
\[
L(1)^{\ominus_\bfi}\cong C_{-\Lm_1}\circ L(2).
\]
For $\bfi'=212$, $L(1)$ corresponds to $(0,1,0)\in \bbB_{212}$. 
One has ${\bf h}_1=(1,1,0)$, ${\bf h}_2=(0,1,1)$. Then, one gets
$(0,-1,0)=\til f_2(0,-1,-1)=\til f_2(-{\bf h}_2)=-{\bf h}_2+\til f_2(0,0,0)$ and then
\[
L(1)^{\ominus_{\bfi'}}\cong C_{-\Lm_2}\circ L(2).
\]
\end{ex}
 
\nd{\bf Problem 2.}
Describe $L\oplus_\bfi L$ and $L^{\ominus_\bfi}$ for 
$L,\,L'\in\tilrgmod$ explicitly.
%What is the corresponding object in $\bbB(\tilrgmod)$ to 
%the summation in $\bbB_\bfi\cong \bbZ^N$?
%That is, 
%for $b,b'\in \bbB_\bfi$, 
%\[
%b+b'\in \bbB_\bfi \longleftrightarrow \wtil\Psi^{-1}(b+b')??\in\bbB(\tilrgmod)
%\]
Indeed, for  ${\bf h}_\Lm\in \cH_\bfi$ and $b\in B(\ify)$
\[
{\bf h}_\Lm+b\in B^{{\bf h}_\Lm}(\ify)\subset \bbB_\bfi \longleftrightarrow C_\Lm\circ \Psi^{-1}(b)\in\bbB(\tilrgmod).
\]
But, note that if $b$ corresponds to "imaginary"(=non-real) simple module 
$S_b$ in $\rgmod$, then $S_b\circ S_b$ is not simple and then 
it does not appear in $\bbB_\bfi$. 

%Fundamental domain $\bbB_\bfi/\cH_\bfi$ ?\\
%\[
%{\bf h}_i\longleftrightarrow \hbox{ lowest weight vector }u_{w_0\Lm_i}
%\in B(\Lm_i)
%\]

\medskip
\subsection{Category $\wtil{\scC}_w$}
In \cite{Loc},  it has been shown that 
for an arbitrary symmetrizable Kac-Moody Lie algebra and 
any Weyl group element  $w\in W$, there exists a subcategory 
$ {\scC}_w\subset$$R$-gmod and it 
admits a localization 
\[
\wtil{\scC}_w=\scC_w[C_i^{\circ-1}\,|\,i\in I],
\q(C_i=M(w\Lm_i,\Lm_i)).
\]
Indeed, note that for finite type Lie algebra setting, ${\scC}_{w_0}=R$-gmod.

\nd {\bf Problem 3.} We conjecture that  the localization $\wtil{\scC}_w$
possess a crystal 
\[
\bbB(\wtil{\scC}_w)=\{S\mid S\hbox{ is a self-dual simple
module in }\wtil\scC_w\}.
\]
 If so, we also conjecture that 
 there is an   isomorphism of crystals
\[
\bbB(\wtil{\scC}_w)\,\mapright{\sim}\,B_{i_1}\ot\cd\ot B_{i_m},
\]
where $i_1\cd i_m$ is a reduced word of $w$.

% {Problem 4}
%The determinantial module $M(w\Lm,u\Lm)\in R$-gmod corresponds to 
%the generalized minor $\Delta_{w\Lm,u\Lm}$. \\
%Can we describe 
%\[
%\wtil\Psi(M(w\Lm,u\Lm))\in \bbB_\bfi
%%\{C_\Lm\,|\,\Lm\in P\}\longleftrightarrow \cH_\bfi
%%\subset \bbB_\bfi
%\] 
%explicitly ? Is there any contribution from this to calculation of the
% generalized minors or $i$-trails?
% 
%\medskip

\subsection{Rigidity}
\begin{df}\label{rigid}
Let $X,Y$ be objects in a monoidal category $\cT$, and $\vep:X\ot Y\to 1$
and $\eta:1\to Y\ot X$ morphisms in $\cT$.
We say that a pair $(X,Y)$ is {\it dual pair} or 
$X$ is a {\it left dual} to  $Y$ or  
$Y$ is a {\it right dual} to  $X$ if 
the following compositions are identities:
\[
%\begin{eqnarray*}
X\simeq X\ot 1\,\,\mapright{{\rm id}\ot \eta}\,\,X\ot Y\ot X
\,\,\mapright{\vep\ot {\rm id}}\,\,1\ot X\simeq X,\,\,
Y\simeq 1\ot Y\,\,\mapright{\eta\ot {\rm id}}\,\,Y\ot X\ot Y
\,\,\mapright{{\rm id}\ot \vep}\,\,Y\ot 1\simeq Y
%\end{eqnarray*}
\]
\end{df}
%\vspace{-20pt}
We denote a right dual to $X$ by $\cD(X)$ and  a 
left dual to $X$ by $\cD^{-1}(X)$.
\vspace{-5pt}
\begin{thm}[\cite{Loc}]
For any finite type $R$,
 $\tilrgmod$ is { rigid}, i.e.,
every object in $\tilrgmod$ has left and right duals.
\end{thm}
Note that in \cite{Loc2}, it is shown that for any symmetrizable Kac-Moody setting
the localized category $\wtil\scC_w$ is rigid.

For a category $C$, let us denote the opposite category of $C$ by $C^{\rm op}$ and 
for a monoidal category $T=(T,\ot)$,   the reversed monoidal category 
$T^{\rm rev}=(T,\ot^{\rm rev})$ is defined by $X\ot^{\rm rev}Y:=Y\ot X$ and
$f\ot^{\rm rev} g:=g\ot f$ for any objects $X,Y\in T$ and any morphisms $f,g$.
%For example, $\cD^{-1}(L(i))=M(w_0\Lm_i,s_i\Lm_i)\circ C_i^{-1}$.
By the results in \cite{Loc,Loc2}, we obtain 
\begin{pro}
The left (resp. right) dual $\cD^{-1}$ (resp. $\cD$) defines the following equivalence of 
categories 
\[
\tilrgmod\longleftrightarrow (\tilrgmod^{\rm op})^{\rm rev}.
\]
\end{pro}
Then by this proposition, 
one finds that for a simple module $S$, both $\cD(S)$ and $\cD^{-1}(S)$ 
are simple. Then, we can consider the following problem. 
\medskip

\nd{\bf Problem 4.}
For a simple object $C_\Lm\circ S\in \bbB(\tilrgmod)$, describe the right and 
left duals explicitly:
\vspace{-5pt}
\[
\wtil\Psi({\cD}(C_\Lm\circ S)),\q 
\wtil\Psi({\cD}^{-1}(C_\Lm\circ S))\in \bbB_\bfi.
\]

%%%%%%%%%%%%%%%%%%%%%%%%%%%%%%%%%%%%%%%%%%%%%%%%%%%%%%%%%%%%%%%%%
\bibliographystyle{amsalpha}

\end{document}